\documentclass[11pt]{article}

\usepackage{amsmath, amssymb, amscd}
\usepackage[matrix,arrow]{xy}

\def\eqref#1{(\ref{#1})}
\newcommand{\goth}{\frak}
\newcommand{\g}{{\frak g}}
\newcommand{\arrow}{{\:\longrightarrow\:}}
\newcommand{\Z}{{\Bbb Z}}
\newcommand{\C}{{\Bbb C}}
\newcommand{\R}{{\Bbb R}}

\renewcommand{\H}{{\Bbb H}}
\newcommand{\6}{\partial}
\renewcommand{\1}{\sqrt{-1}\:}

\renewcommand{\c}[1]{{\cal #1}}
\newcommand{\calo}{{\cal O}}

\renewcommand{\tilde}{\widetilde}
\renewcommand{\bar}{\overline}
\renewcommand{\phi}{\varphi}
\renewcommand{\epsilon}{\varepsilon}

\renewcommand{\leq}{\leqslant}

\newcommand{\ev}{{\rm even}}
\newcommand{\odd}{{\rm odd}}

\newcommand{\End}{\operatorname{End}}

\newcommand{\Vol}{\operatorname{Vol}}
\newcommand{\Hom}{\operatorname{Hom}}

\newcommand{\comment}[1]{{}}

\def\blacksquare{\hbox{\vrule width 4pt height 4pt depth 0pt}}
\def\endproof{\blacksquare}

\makeatletter

\@ifundefined{Bbb}
     {\newcommand{\Bbb}[1]{{\mathbb #1}}}%
{}

\newcommand{\ps@verbit}{%
  \renewcommand{\@oddhead}{%
          \scriptsize
          {Hodge theory on HKT-manifolds}
          \hfil\tiny {M. Verbitsky, \ \ \ December 20, 2001 }}
  \renewcommand{\@evenhead}{\@oddhead}
  \renewcommand{\@oddfoot}{\hfil\thepage\hfil}
  \renewcommand{\@evenfoot}{\@oddfoot}}
 
\newcounter{Mycounter}[section]
\newcounter{lemma}[section]
\renewcommand{\thelemma}{{Lemma \thesection.\arabic{lemma}}}
\newcommand{\lemma}{%
     \setcounter{lemma}{\value{Mycounter}}
     \refstepcounter{lemma}
     \stepcounter{Mycounter}
     {\bf \thelemma:\ }}

\newcounter{claim}[section]
\renewcommand{\theclaim}{{Claim \thesection.\arabic{claim}}}
\newcommand{\claim}{%
     \setcounter{claim}{\value{Mycounter}}
     \refstepcounter{claim}
     \stepcounter{Mycounter}
     {\bf \theclaim:\ }}

\newcounter{sublemma}[section]
\newcounter{corollary}[section]
\renewcommand{\thecorollary}{{Corollary \thesection.\arabic{corollary}}}
\newcommand{\corollary}{%
     \setcounter{corollary}{\value{Mycounter}}
     \refstepcounter{corollary}
     \stepcounter{Mycounter}
     {\bf \thecorollary:\ }}

\newcounter{theorem}[section]
\renewcommand{\thetheorem}{{Theorem \thesection.\arabic{theorem}}}
\newcommand{\theorem}{%
     \setcounter{theorem}{\value{Mycounter}}
     \refstepcounter{theorem}
     \stepcounter{Mycounter}
     {\bf \thetheorem:\ }}

\newcounter{conjecture}[section]
\newcounter{proposition}[section]
\renewcommand{\theproposition}
       {{Proposition \thesection.\arabic{proposition}}}
\newcommand{\proposition}{%
     \setcounter{proposition}{\value{Mycounter}}
     \refstepcounter{proposition}
     \stepcounter{Mycounter}
     {\bf \theproposition:\ }}

\newcounter{definition}[section]
\renewcommand{\thedefinition}
       {{Definition~\thesection.\arabic{definition}}}
\newcommand{\definition}{%
     \setcounter{definition}{\value{Mycounter}}
     \refstepcounter{definition}
     \stepcounter{Mycounter}
     {\bf \thedefinition:\ }}

\newcounter{example}[section]
\renewcommand{\theexample}{{Example \thesection.\arabic{example}}}
\newcommand{\example}{%
     \setcounter{example}{\value{Mycounter}}
     \refstepcounter{example}
     \stepcounter{Mycounter}
     {\bf \theexample:\ }}

\newcounter{remark}[section]
\renewcommand{\theremark}{{Remark \thesection.\arabic{remark}}}
\newcommand{\remark}{%
     \setcounter{remark}{\value{Mycounter}}
     \refstepcounter{remark}
     \stepcounter{Mycounter}
     {\bf \theremark:\ }}

\newcounter{problem}[section]
\newcounter{question}[section]
\@addtoreset{equation}{section}
\@addtoreset{footnote}{section}
\makeatother

\begin{document}

\begin{center}
{\LARGE\bf
Hyperk\"ahler manifolds with torsion,\\[3mm] supersymmetry and Hodge theory
}
\\[4mm]
Misha Verbitsky,\footnote{The author is 
partially supported by CRDF grant RM1-2087}
\\[4mm]
{\tt verbit@dnttm.rssi.ru, \ \  verbit@mccme.ru}
\end{center}

{\small 
\hspace{0.15\linewidth}
\begin{minipage}[t]{0.7\linewidth}
{\bf Abstract} \\
Let $M$ be a hypercomplex Hermitian manifold,
$(M,I)$ the same manifold considered as a complex
Hermitian with a complex structure $I$ induced by
the quaternions. The standard linear-algebraic construction produces
a canonical nowhere degenerate (2,0)-form $\Omega$ 
on $(M,I)$. It is well known that $M$ is hyperk\"ahler
if and only if the form $\Omega$ is closed. 
The $M$ is called HKT (hyperk\"ahler with torsion)
if $\Omega$ is closed with respect to the Dolbeault
differential $\6:\; \Lambda^{2,0}(M,I) \arrow \Lambda^{3,0}(M,I)$.
Conjecturally, all compact hypercomplex manifolds
admit an HKT-metrics. We exploit a remarkable 
analogy between the de Rham DG-algebra 
of a K\"ahler manifold and the Dolbeault DG-algebra of an HKT-manifold.
The supersymmetry of a K\"ahler manifold $X$ is given by 
an action of an 8-dimensional Lie superalgebra $\g$ on
$\Lambda^*(X)$, containing the Lefschetz $SL(2)$-triple,
the Laplacian and the de Rham differential. We establish
the action of $\g$ on the Dolbeault DG-algebra
$\Lambda^{*,0}(M,I)$ of an HKT-manifold. This
is used to construct a canonical Lefschetz-type
$SL(2)$-action on the space of harmonic spinors of $M$.
\end{minipage}
}

{
\small
\tableofcontents
}

\section{Introduction}
\label{_Intro_Section_}


Hyperk\"ahler manifolds with torsion (HKT-manifolds) were
introduced by P.S.Howe and G.Papadopoulos (\cite{_Howe_Papado_})
were much studied in physics literature since then.
For an excellent survey of these works written from  a mathematician's
point of view, the reader is referred to the paper of
G. Grantcharov and Y. S. Poon \cite{_Gra_Poon_}.
In physics, HKT-manifolds appear as moduli of
brane solitons in supergravity and M-theory (\cite{_GP2_}, 
\cite{_P:Rome_}). HKT-manifolds also arise as moduli
space of some special black holes in N=2 supergravity 
(\cite{_GP1_}, \cite{_GPS_}).

The term ``hyperk\"ahler manifold with torsion'' is actually
quite misleading, because an HKT-manifold is not hyperk\"ahler.
This is why we prefer to use the abbreviation ``HKT-manifold''.

\subsection{Hypercomplex manifolds}

Let $M$ be a smooth manifold
equipped with an action of 
quaternion algebra on its tangent bundle.
The manifold $M$ is called {\bf hypercomplex} (\cite{_Boyer_})
if for any quaternion $L\in \H$, $L^2=-1$,
the corresponding almost complex structure is 
integrable. If, in addition, $M$ admits a 
Riemannian structure, and 
for any quaternion $L\in \H$, $L^2=-1$,
$L$ establishes a K\"ahler structure on $M$,
the manifold $M$ is called {\bf hyperk\"ahler}.

The geometry of hypercomplex manifold
is quite rich, but not completely understood yet.
There are compact homogeneous examples
(\cite{_Joyce_}), compact
inhomogeneous examples (\cite{_Poon_Pedersen_}).
The task of producing new examples of compact hypercomplex
manifolds is clearly not as difficult as 
that for compact hyperk\"ahler manifolds.

The main tool of hypercomplex geometry is the so-called
{\bf Obata connection}. Given a hypercomplex manifold $M$,
there exists a unique torsion-free connection on the
tangent bundle $TM$ which preserves the quaternion action;
this connection was introduced by M. Obata in 1950-ies
(\cite{_Obata_}). Clearly, the holonomy group
of the Obata connection lies in $GL(n, \H)$,
where $n = \dim_\H M$. The manifold $M$ admits
a hyperk\"ahler metrics if and only if its 
holonomy preserves a positive definite
metrics. 

There are many results on deformations and 
Dolbeault cohomology
of a homogeneous hypercomplex manifold
(\cite{_Poon_Pedersen_}, \cite{_GPP_}, etc.)
However, the algebraic geometry of a general
hypercomplex manifold is {\it terra incognita},
so far. 

Since these manifolds are (usually) not K\"ahler, no
relation between the de Rham and Dolbeault cohomology is
established. This is why the most natural 
geometrical questions are so difficult to solve.

\subsection{HKT-metrics on hypercomplex manifolds}

For a historically correct definition of an HKT-manifold
(\cite{_Howe_Papado_}), see
Subsection \ref{_HKT_defi_Subsection_}. However, throughout
this paper, we use not this definition, but its reformulation,
introduced by Grantcharov and Poon (\cite{_Gra_Poon_}).

\hfill

Let $M$ be a hypercomplex manifold, and $h$ a Riemannian
metrics on $M$. The metrics $h$ is called {\bf quaternionic Hermitian}
if $h$ is invariant with respect to the group
$SU(2)\subset \Bbb H^*$ of unitary quaternions.
In this case, $M$ is called {\bf a hypercomplex
Hermitian manifold} (\ref{_Hermitian_hype_Definition_}). 

Given a hypercomplex Hermitian manifold, and an
induced complex structure $L\in \Bbb H$, $L^2=-1$,
one may consider the corresponding real-valued $(1,1)$-form
$\omega_L \in \Lambda^{1,1}(M, L)$,
$\omega_L(x,y) := h(x, L y)$. This is a well
known anti-symmetric form associated with the
Hermitian structure on a complex manifold $(M, L)$.

Let $I, J, K\in \Bbb H$ be the standard basis in quaternions.
Consider the differential form $\Omega:= \frac 1 2 (\omega_J + \1 \omega_K)$.
An elementary linear-algebraic calculation ensures that $\Omega$
is of type $(2,0)$ with respect to the complex structure $I$:
\[ \Omega \in \Lambda^{2,0}(M,I).\]

The following lemma
 is due to Grantcharov and Poon
(it is a reformulation of Proposition 2 of
\cite{_Gra_Poon_}).

\hfill

\lemma \label{_HKT_cond_Lemma_}
Let $M$ be a hypercomplex Hermitian manifold.
Consider the subspace \[ V\subset \Lambda^3(M)\]  generated by
3-forms of type \[ d\omega_L, \ \ L\in \Bbb H, \ \ L^2=-1.\] 
The space $\Lambda^3(M)$ is equipped with 
a natural action of the  group of unitary quaternions 
$SU(2)\subset {\Bbb H}^*$. Let  $I, J, K\in \Bbb H$ be the
standard triple of quaternions.
Then the following conditions are equivalent.
\begin{description}
\item[(i)] The space $V$ belongs to a direct sum of several
irreducible $SU(2)$-subresentations $V_i\subset \Lambda^3(M)$
of dimension 2. \footnote{In other words, all elements
of $V$ have weight 1 with respect to the $SU(2)$-action.}

\item[(ii)] We have $\6\Omega=0$, where
$\Omega\in \Lambda^{2,0}(M,I)$, 
$\Omega:= \frac 1 2 (\omega_J + \1\omega_K)$
is the $(2,0)$-form constructed above,
and $\6:\; \Lambda^{2,0}(M,I)\arrow \Lambda^{3,0}(M,I)$
the Dolbeault differential. 
\end{description}

{\bf Proof:} \ref{_HKT_via_diff_of_Omega_Theorem_} (see also
\cite{_Gra_Poon_}, Proposition 2). \endproof

\hfill

\definition
Let $M$ be a hypercomplex Hermitian manifold, 
$I, J, K\in \Bbb H$ the standard basis in quaternions,
and \[ \Omega\in \Lambda^{2,0}(M,I),\ \  \Omega:= \frac 1 2 (\omega_J + \1\omega_K)\]
the $(2,0)$-form constructed above.
Then $M$ is called {\bf HKT-manifold} (or {\bf hyperk\"ahler
with torsion}) if
\begin{equation}\label{_HKT_intro_Equation_}
\6\Omega=0.
\end{equation}

\hfill

\remark 
By \ref{_HKT_cond_Lemma_}, the HKT-condition does not
depend from the choice of a quaternion
basis $I, J, K\in \Bbb H$. 

\hfill

The form $\Omega$ is clearly non-degenerate. The condition
\eqref{_HKT_intro_Equation_} means that $\Omega$ is
``Dolbeault symplectic''. This paper is an attempt to 
exploit the analogy between symplectic geometry and 
the HKT-geometry.

\hfill

So far, almost every constructed example of a compact hypercomplex
manifold comes bundled with a natural HKT-metrics.
Moreover, there are no examples of a
compact hypercomplex manifold where
non-existence of an HKT-metrics is established. Therefore,
HKT-manifolds are a natural object to study.

\hfill

The most prominent example of an HKT-manifold 
is due to D. Joyce (\cite{_Joyce_}) and 
Spindel et al (\cite{_SSTvP_})
who constructed the hypercomplex structures
on a compact Lie group, and Opfermann-Papadopoulos, 
who noticed that this manifold admits an HKT-metrics
  (see \cite{_OP_}, \cite{_Gra_Poon_}). 
Let $G$ be a compact semisimple Lie group.
Then there exists a number $n$, depending 
on the geometry of a group $G$, such that
the product of $G$ and an $n$-dimensional compact torus
admits a left-invariant hypercomplex structure 
(\ref{_Joyce_exa_Example_}).
The HKT-metrics is provided by the Killing form.

\subsection{The K\"ahler-de Rham superalgebra}

In this paper, we study the Dolbeault complex of an HKT-manifold,
from the viewpoint of Hodge theory and Kodaira relations.
In K\"ahler geometry, the Kodaira relations might be interpreted as
relations in a certain Lie superalgebra (see
Subsection \ref{_LSA_Subsection_} for a definition
of the Lie superalgebras). To be more precise,
let $X$ be a K\"ahler manifold, $L, \Lambda, H$ the
Lefschetz $SL(2)$-triple (\cite{_Griffi_Harri_}),
and $\6, \bar\6$ the Dolbeault differentials,
considered as odd endomorphisms of the graded 
vector space $\Lambda^*(M)$.
Consider the Lie superalgebra $\g \subset \End(\Lambda^*(M))$
generated by $L, \Lambda, H$ (even), $\6$ and $\bar\6$ (odd).
 The Kodaira relations are interpreted as the relations
in the Lie superalgebra $\g$ 
(see Section \ref{_KdR_Section_}). From these relations
it follows that the Lie algebra $\g$ is independent from
the choice of a K\"ahler manifold $X$.
We call $\g$ {\bf the K\"ahler-de Rham  superalgebra}.

\subsection{The K\"ahler-de Rham  superalgebra and HKT-geometry}

Let $M$ be an HKT-manifold, $I$, $J$, $K\in \Bbb H$
the standard triple of induced complex structures,
and $(\Lambda^{*,0}(M, I), \6)$ the Dolbeault
differential graded (DG-) algebra of $(M,I)$. 
We use an analogy between 
$(\Lambda^{*,0}(M, I), \6)$
and the de Rham DG-algebra
of a K\"ahler manifold. The role of the de Rham
differential is played by the Dolbeault differential
on $\Lambda^{*,0}(M, I)$. The role of the K\"ahler
form $\omega\in \Lambda^{1,1}(X)$ is played 
by the $(2,0)$-form $\Omega\in \Lambda^{2,0}(M, I)$
constructed above. One can associate with
$\Omega$ an $SL(2)$-triple $L_\Omega, \Lambda_\Omega,
H_\Omega$ acting on $\Lambda^{*,0}(M, I)$ 
(\ref{_K-dR_inside_hK-dR_Dolbeault_Proposition_}).

We have the complex structure operator
$I:\; \Lambda^*(X) \arrow \Lambda^*(X)$
acting on the differential forms over a  
K\"ahler manifold $X$ in a multiplicative way.
On the hypercomplex side, its role is played
by an operator $\c J$ defined as follows.

Let $J:\; \Lambda^*(M) \arrow \Lambda^*(M)$
be the complex structure operator
induced by $J\in \Bbb H$. Since 
the quaternions $I$ and $J$ anti-commute,
$J$ maps $\Lambda^{p,q}(M, I)$ to $\Lambda^{q,p}(M, I)$.
Composing $J$ with the complex conjugation,
we obtain an anticomplex automorphism
\[ \c J:\; \Lambda^{p,q}(M, I) \arrow \Lambda^{p,q}(M, I).
\]
We consider $\c J$ as an endomorphism
of $\Lambda^{*,0}(M, I)$.

Twisting the Dolbeault differential with $\c J$, we obtain
a differential $\6_J$, which is an analogue of the
twisted de Rham differential $d^c$ on a K\"ahler manifold.
Consider the Lie superalgebra generated by the 
Lefschetz triple $L_\Omega, \Lambda_\Omega,
H_\Omega$ and the differentials $\6, \6_J$.
We show that this  Lie superalgebra is
isomorphic to the K\"ahler-de Rham superalgebra
(\ref{_K-dR_on_HKT_Corollary_}). 
In fact, this isomorphism is equivalent to the
HKT-condition.

When the manifold $M$ is not only HKT, but also hyperk\"ahler,
the Dolbeault Laplacian $\6\6^*+ \6^* \6$ is a central
element of the K\"ahler-de Rham superalgebra constructed
above. Therefore, we have a Lefschetz-type $SL(2)$-action
on the Dolbeault cohomology $H^{*,0}(M,I)$ of $(M,I)$. This is
not very useful by itself, because $H^{*,0}(M,I)$ 
of a compact hyperk\"ahler manifold is very easy
to compute. However, the same superalgebra
acts on differential forms with coefficients
in appropriate vector bundles (see 
\cite{_Verbitsky:Hyperholo_bundles_}).
This is used to construct an $SL(2)$-action
on the cohomology of these vector bundles.
As an end result, we show that 
the deformation space of hyperholomorphic
vector bundles is unobstructed, and 
prove that it is a hyperk\"ahler variety
(see \cite{_Verbitsky:Hyperholo_bundles_}, \cite{_Verbitsky:hypercomple_},
\cite{_KV:book_} for details).

A similar approach cannot work in whole generality for
an HKT-\-ma\-ni\-fold. Indeed, if we have a Lefschetz-type
action on the Dolbeault cohomology of $(M,I)$,
we have an isomorphism $H^0(\calo_{(M,I)})\cong H^n(\calo_{(M,I)})$,
where $n=\dim_\C M$. By Serre's duality,
$H^n(\calo_{(M,I)})\cong H^0({\sf K})^*$,
where ${\sf K}:=K_{(M,I)}$ is the canonical bundle of $(M,I)$.
This bundle is known to be non-trivial in some examples
(for the Hopf surface, for example), and to have no
non-trivial sections. But then,
the 1-dimensional space $H^0(\calo_{(M,I)})$ would
be isomorphic to the trivial space $H^0({\sf K})^*$,
which is impossible.
Therefore, one cannot expect an analogue of 
Lefschetz theorem to hold for
the Dolbeault cohomology of $(M,I)$.

However, it is natural to expect some kind of Hodge 
theory to hold on the space 
$\Lambda^{*,0}(M,I)\otimes {\sf K}^{1/2}$
of ``$(p,0)$ half-forms''. First, by Serre's duality, we have
\[ H^i({\sf K}^{1/2}) \cong H^{n-i}({\sf K}^{1/2})^*,
\]
so the Lefschetz action identifies the spaces of the
same dimension. Second, 
the space $\Lambda^{*,0}(M, I) \otimes {\sf K}^{1/2}$
is naturally identified with the space of spinors on $M$,
which is a Riemannian invariant of $M$. The corresponding
cohomology space $H^*({\sf K}^{1/2})$ is the space 
of harmonic spinors, and it is also a Riemannian invariant.
Thus, the half-forms and the cohomology of
${\sf K}^{1/2}$ are in some sense more relevant
geometrically than the Dolbeault cohomology.

This paper is dedicated to revealing the
Lefschetz-type action on the harmonic spinors.

\hfill

First of all, we notice that the canonical bundle
$\sf K$ of $(M,I)$ is topologically trivial. Indeed, taking
the top exterior power of the symplectic form $\Omega\in \Lambda^{2,0}(M,I)$,
we obtain a nowhere degenerate section of $\sf K$. Using
this trivialization, one defines the square root
of $\sf K$ in a usual way (see 
Subsection \ref{_Lefschetz_Subsection_}).
The endomorphism 
\[ \c J:\; \Lambda^{p,0}(M,I)\otimes {\sf K}^{1/2} \arrow 
\Lambda^{p,0}(M,I)\otimes {\sf K}^{1/2}
\]
is defined in the same fashion as above.
Let ${}^n\6$ be the Dolbeault differential in
$\Lambda^{*,0}(M,I)\otimes {\sf K}^{1/2}$,
\[ {}^n\6:\; \Lambda^{p,0}(M,I)\otimes {\sf K}^{1/2}\arrow
   \Lambda^{p+1,0}(M,I)\otimes {\sf K}^{1/2},
\]
and ${}^n\6_J$ the twisted differential, 
\[ {}^n\6_J := - \c J\circ {}^n\6\circ \c J.
\]
The Lefschetz $SL(2)$-triple $L_\Omega$, $\Lambda_\Omega$,
$H_\Omega$ acts on $\Lambda^{*,0}(M,I)\otimes {\sf K}^{1/2}$
in a natural way. The main result of this paper
is the following

\hfill

\theorem
 [ Subsection \ref{_normalized_g_defi_Subsection_} ]
Let $M$ be a HKT-manifold, $I,J,K\in \Bbb H$ the standard
basis in quaternions, and ${\sf K}^{1/2}$ the square
root of the canonical bundle of $(M,I)$, constructed above.
Consider the Lefschetz $SL(2)$-triple $L_\Omega$, $\Lambda_\Omega$,
$H_\Omega$ acting on the bundle of $(p,0)$-half-forms 
$\Lambda^{*,0}(M,I)\otimes {\sf K}^{1/2}$ defined above.
Let
\[ {}^n\6:\; \Lambda^{p,0}(M,I)\otimes {\sf K}^{1/2}\arrow
   \Lambda^{p+1,0}(M,I)\otimes {\sf K}^{1/2},
\]
be the Dolbeault operator, and ${}^n\6_J$ the twisted differential, 
\[ {}^n\6_J := - \c J\circ {}^n\6\circ \c J.
\]
Consider the Lie superalgebra $\g$ generated by the 
even operators $L_\Omega$, $\Lambda_\Omega$,
$H_\Omega$ and the odd operators ${}^n\6$, ${}^n\6_J$,
Then $\g$ is isomorphic to the K\"ahler-de Rham superalgebra.

\hfill

{\bf Proof:} See Subsection \ref{_normalized_g_defi_Subsection_}.
\endproof

\hfill

\corollary 
We have a Lefschetz-type $SL(2)$-action on the cohomology
$H^*({\sf K}^{1/2})$. 

\hfill

{\bf Proof:} This is \ref{_Lefschetz_Theorem_}. \endproof

\hfill

This result should be especially useful when the canonical bundle
${\sf K}$ of $(M,I)$ is trivial. However, we have not found 
examples of compact HKT- (or even hypercomplex)
manifolds which are not hyperk\"ahler
and have trivial canonical bundle.
Moreover, in all non-hyperk\"ahler
examples where we have computed the
group $H^i({\sf K}^{1/2})$, this group was trivial,
for all $i$. Still, the action of the K\"ahler-de Rham superalgebra
on the space of spinors is quite remarkable.

\subsection{Non-K\"ahler manifolds and Calabi-Yau geometry}

As we have seen above, the canonical bundle of a hypercomplex
manifold is topologically trivial. Therefore, these manifolds
may be regarded as ``non-K\"ahler Calabi-Yau manifolds''.
The study of such manifolds is quite important,
due to the following conjecture of M. Reid. 
Let $X$ be a 3-dimensional Calabi-Yau manifold.
Suppose that $X$ contains a rational curve which can be
blown down; denote the blow-down of $X$ by $\tilde X_1$.
Consider a generic small deformation $X_1$ of $\tilde X_1$. 
It is known that all singularities of 3-dimensional
Calabi-Yau manifolds can be removed by a 
small deformation. Therefore, we may assume
that $X_1$ is a smooth manifold, which has trivial canonical class
and satisfies $h^2(X_1) = h^2(X)-1$. Repeating this process,
we obtain a manifold $X_n$ with $h^2(X_n)=0$. 
This manifold is, of course, non-K\"ahler. 
However, its topological structure is very simple. 
Namely, $X_n$ is diffeomorphic to a connected sum 
$\#_k S^3\times S^3$ of $k$ copies of 
a $S^3\times S^3$, where $S^3$ is a 
3-dimensional sphere, and $k = h^3(X)/2$.

M. Reid conjectured that, starting from
another 3-dimensional Calabi-Yau manifold
$Y$, $h^3(Y)=h^3(X)$, we obtain a non-K\"ahler manifold
$Y_{n'}$ in the same deformation class as $X_n$,

This conjecture is significant because, in Calabi-Yau geometry,
most results deal with
complete intersections in toric and homogeneous manifolds,
and the general Calabi-Yau manifolds are quite 
difficult to study. The M. Reid's conjecture
gives a possibility to reduce a given Calabi-Yau manifold to a 
toric Calabi-Yau  by a series of birational
transforms and deformations. This way, one
might hope to extend the
standard results about 
complete intersections in toric manifolds
(the Mirror Conjecture, for instance)
to the general case.

The M. Reid's conjecture is difficult to tackle because
very little is known about the geometry of non-K\"ahler
Calabi-Yau manifolds. The present paper can be read as an attempt
to study their geometry, from the Hodge-theoretic point of view.
The analogy is straightforward. The main
working example of an HKT manifold, the 
compact group $SU(3)$, looks very
similar, from the topological 
and geometrical point of view,
to the manifold $\#_k S^3\times S^3$.
The simplest 3-dimensional
``non-K\"ahler Calabi-Yau''
manifold $S^3\times S^3$ is actually a 
Lie group; and the Lie group $S^3\times S^1$ is
also an HKT-manifold, called the Hopf surface.

\subsection{Contents}

\begin{itemize}
\item This Introduction (Section \ref{_Intro_Section_})
is independent from the rest of this paper.

\item In Section \ref{_hyperka_Section_},
we recite some basic results and conventions of hypercomplex
and hyperkaehler geometry. This material is standard
(see, e.g., \cite{_KV:book_}).

\item In Section \ref{_HKT_Section_},
we repeat the definition and some basic
properties of HKT-manifolds (\cite{_Howe_Papado_}, 
\cite{_Gra_Poon_}). 

\item In Section \ref{_KdR_Section_}, 
we study the natural superalgebra of a K\"ahler manifold
(\cite{_FKS_}).

\item In Section \ref{_q.D._Section_},
we define and study the quaternionic
Dolbeault complex, following S. Salamon
(see \cite{_Capria-Salamon_}
and \cite{_V:projective_}). 

\item In Section \ref{_hk-DR_Section_}, 
we study the supersymmetry of a hyperk\"ahler manifold,
following \cite{_FKS_} and \cite{_Verbitsky:SO(5)_}.

\item The material of Sections 
\ref{_hyperka_Section_}-\ref{_hk-DR_Section_} is
known from the literature; results of Sections 
\ref{_superalh_HKT_Section_}-\ref{_K^1/2_Section_} are new.

\item In Section \ref{_superalh_HKT_Section_},
we prove that the Dolbeault algebra of an HKT-manifold
admits the same kind of supersymmetry as the de Rham algebra
of a K\"ahler manifold. 

\item In Sections \ref{_Kodaira_Section_}-\ref{_*_Section_}, we 
obtain some commutation relations in the superalgebra
of the Dolbeault complex of an HKT-manifold. 

\item In Section \ref{_K^1/2_Section_}, 
we apply these commutation relation to establish the geometrically
relevant kind of supersymmetry on the bundle of spinors. This
is used to show that the harmonic spinors admit 
a natural Lef\-schetz-\-type $SL(2)$-action

\end{itemize}


\section{Hypercomplex manifolds}
\label{_hyperka_Section_}


\subsection{Hypercomplex Hermitian manifolds}

\definition \label{_hc_Definition_}
Let $M$ be a smooth manifold equipped with an action
of the quaternion algebra $\Bbb H$ in its tangent space.
Assume that for all quaternions $L\in \Bbb H$, $L^2=-1$,
the almost complex structure given by $L\in \End(TM)$
is integrable. Then $M$ is called {\bf hypercomplex}.

\hfill

\remark
By D. Kaledin's theorem
(\cite{_Kaledin_}), in \ref{_hc_Definition_} it
suffices to check integrability only for
$L_1, L_2\in \Bbb H$, $L_1\neq \pm L_2$.
The integrability of $L_1$, $L_2$
 implies integrability of any
almost complex structure given by
$L\in \Bbb H$, $L^2=-1$.

\hfill

This paper is dedicated to the study of hypercomplex
manifolds. For this purpose, we introduce a natural
kind of Riemannian metrics, called quaternionic
Hermitian. It is related to the hypercomplex
structure in the same way as the usual Hermitian
metrics on a complex manifold is related to a
complex structure. This definition is purely 
linear-algebraic.

\hfill

Let $V$ be a quaternionic vector space. Given $L\in \Bbb H$,
$L^2=-1$, $L$ defines a complex structure on $V$. We denote
$V$, considered as a complex vector space, by $(V,L)$. The following
lemma is trivial. 

\hfill

\lemma \label{_q_Hermi_Lemma_}
Let $V$ be a quaternionic vector space, and $h:\; V \times V \arrow \R$
a positive bilinear symmetric form. Then the following conditions
are equivalent.

\begin{description}
\item[(i)] For any $L\in \Bbb H$,
$L^2=-1$, the metrics $h$ is Hermitian on the complex vector space 
$(V,L)$.
\item[(ii)] Consider the action of the group $SU(2)$ of unitary
quaternions on $V$. Then $h$ is $SU(2)$-invariant.
\end{description}

{\bf Proof:} Clear. \endproof

\hfill

\definition
Let $V$ be a quaternionic vector space, and $h:\; V \times V \arrow \R$
a positive bilinear symmetric form. Then space $(V, h)$ is called 
{\bf a quaternionic Hermitian space} if any of the conditions
of \ref{_q_Hermi_Lemma_} holds. In this case, the metrics $h$
is also called {\bf quaternionic Hermitian}.

\hfill

\definition \label{_Hermitian_hype_Definition_}
Let $M$ be a hypercomplex manifold, and $s$ a Riemannian structure
on $M$. We say that $M$ is {\bf hypercomplex Hermitian} if
for all $x\in M$, the tangent space $(T_xM, s)$ is quaternionic 
Hermitian.

\hfill

\definition \label{_induced_c_str_Definition_}
Let $M$ be a hypercomplex manifold, $L\in \Bbb H$,
$L^2=-1$. By definition, $L$ induces a complex structure on
$M$. Thus obtained complex manifold is denoted by $(M,L)$.
The complex structure $L$ is called {\bf induced by
the hypercomplex structure}. If, in addition,
$M$ is hypercomplex Hermitian, then $(M,L)$
is a complex Hermitian manifold. The 
corresponding real-valued skew-symmetric 2-form
is denoted by $\omega_L\in \Lambda^{1,1}(M,L)$,
\[ \omega_L(x,y) = s(x, Ly),
\]
where $s$ denotes the Riemannian form.

\hfill

\definition
Let $M$ be a hypercomplex Hermitian manifold.
Then $M$ is called {\bf hyperk\"ahler} if 
for any induced complex structure $L\in \Bbb H$,
$L^2=-1$, the corresponding Hermitian manifold
$(M,L)$ is K\"ahler. 

\hfill

Let $M$ be a hypercomplex Hermitian manifold,
and $I,J, K\in \Bbb H$ the standard triple of quaternions.
Consider the 2-form
\begin{equation}\label{_sta_2,0_form_Equation_}
\Omega:=\; \frac 1 2 (\omega_J + \1\omega_K). 
\end{equation}
An elementary linear-algebraic computation shows that
$\Omega$ is of type $(2,0)$ with over $(M,I)$.

\hfill

\claim
Let $M$ be a hypercomplex Hermitian manifold,
$I,J, K\in \Bbb H$ the standard triple of quaternions,
and $\Omega\in \Lambda^{2,0}(M,I)$ the (2,0)-form
constructed above. Then $M$ is hyperk\"ahler if
and only if $d\Omega=0$.

\hfill

{\bf Proof:} The proof is well known
(see \cite{_Besse:Einst_Manifo_}). \endproof

\subsection{The natural $SU(2)$-action on the differential forms}

Let $M$ be a hypercomplex manifold. We identify the group $SU(2)$
with the group of unitary quaternions. This gives a canonical 
action of $SU(2)$ on the tangent bundle, and all its tensor
powers. In particular, we obtain a natural action of $SU(2)$
on the bundle of differential forms. 

The corresponding Lie algebra action is related to the usual Hodge 
decomposition as follows. Let $L$ be an induced complex structure,
and
\[ \Lambda^*(M) = \oplus \Lambda^{p,q}(M,L)
\]
be the associated Hodge decomposition. Consider an operator
\[ ad L:\; \Lambda^i(M) \arrow \Lambda^i(M)\] acting
as $\eta \arrow (p-q)\1 \eta$ on $\eta \in \Lambda^{p,q}(M,L)$.

\hfill

\lemma \label{_su(2)_action_explici_Lemma_}
Let $M$ be a hypercomplex manifold,
$I,J, K\in \Bbb H$ the standard triple of quaternions,
and $ad I, \ ad J, \  ad K$ the corresponding
endomorphisms of $\Lambda^*(M)$. Then 
$ad I, \ ad J, \  ad K$ generate 
the Lie algebra $\goth{su}(2)\subset \End(\Lambda^*(M))$
associated with the hypercomplex structure.

\hfill

{\bf Proof:} Clear (see e.g. 
\cite{_Verbitsky:SO(5)_}, \cite{_Verbitsky:Hyperholo_bundles_}). \endproof


\section{HKT-manifolds}
\label{_HKT_Section_}


\subsection{Busmut connections}

For a reference and a bibliography on Bismut connections,
see \cite{_Gra_Poon_}.

Let $M$ be a complex manifold equipped with a Hermitian metrics, and 
$\nabla$ a connection (not necessarily torsion-free)
on $TM$ which preserves the Hermitian metrics and the complex 
structure. 
Denote its torsion by $T_\nabla \in \Lambda^2 T^* M \otimes TM$.
Using the Riemannian structure to identify 
$TM$ and $T^*M$, we may consider $T_\nabla$ 
as an element of 
$\Lambda^2M \otimes \Lambda^1M\subset (T^*M)^{\otimes 3}$.

\hfill

\definition \label{_Bism_conne_Definition_}
The connection $\nabla$ is called {\bf Bismut connection}
if the tensor 
\[ T_\nabla\in \Lambda^2M 
\otimes \Lambda^1M\subset (T^*M)^{\otimes 3}
\]
is totally skew-symmetric, that is, belongs to the
space of differential forms 
$\Lambda^3(M)\subset (T^*M)^{\otimes 3}$.

\hfill

\theorem\label{_Bismut_conne_exi_Theorem_}
(Chern)
Let $M$ be a complex Hermitian manifold. 
Then the Bismut connection exists and is unique.

{\bf Proof:} Well known (\cite{_Gra_Poon_}). The proof
is constructed along the same lines
as the proof of the existence and uniqueness
of a Levi-Civita connection on a Riemannian manifold.
The space of connections is affine;
one considers the map $\nabla\arrow T_\nabla$
as an affine map of affine spaces, and 
estimates its kernel and cokernel by
dimension count. \endproof

\hfill

\remark 
Let $M$ be a complex Hermitian manifold, $\nabla$ the Bismut connection,
and $T_\nabla \in \Lambda^3(M)$ the torsion tensor, considered as
a differential form as in \ref{_Bism_conne_Definition_}.
It is possible to express $T_\nabla$ in terms of the standard
skew-symmetric form $\omega\in \Lambda^{1,1}(M)$ associated
with the Hermitian form. Namely, let 
$I:\; \Lambda^3(M) \arrow \Lambda^3(M)$ be the
complex structure operator extended to
$\Lambda^3(M)$ by multiplicativity. Then
$T_\nabla = I(d\omega)$ (\cite{_Gra_Poon_}).

\subsection{HKT-manifolds: the definition}
\label{_HKT_defi_Subsection_}

HKT-manifolds were introduced by P.S.Howe and G.Papadopoulos
(\cite{_Howe_Papado_}).
For a reference and a bibliography on HKT-manifolds,
see \cite{_Gra_Poon_}.

\hfill

\definition\label{_HKT_Definition_}
\cite{_Howe_Papado_} Let $M$ be a hypercomplex Hermitian manifold 
(\ref{_Hermitian_hype_Definition_}).
For an induced complex structure $L = I, J, K$, consider
the corresponding 
complex Hermitian manifold $(M,L)$.
Let $\nabla_L$ be the associated Bismut connection
(\ref{_Bismut_conne_exi_Theorem_}). 
We say that $M$ is {\bf HKT- (hyperk\"ahler with torsion)
manifold}  if
\[ \nabla_I = \nabla_J = \nabla_K. \]

\hfill

\remark
An HKT-manifold is {\it not} hyperk\"ahler. Therefore,
the term ``hyperk\"ahler with torsion'' is not legitimate
and actually misleading.
This is why we, throughout this paper, prefer
to use the term ``HKT-manifold'' instead.

\hfill

\example\label{_Joyce_exa_Example_}
Let $G$ be a compact Lie group from the following list:
\begin{equation}\label{_Joyce_list_Equation_}
\begin{aligned}
{} &SU(2l+1), T^1\times SU(2l), \ \ \ T^l \times SO(2l+1),\ \ \  T^l \times Sp(l),\\
{} &  T^{2l-1}\times SO(4l+2),\ \ \  T^2 \times E_6, \ \ \ T^7 \times E_7, \\
{} &T^4 \times F_4,\ \ \  T^2\times G_2.
\end{aligned}
\end{equation}
where $T^i= (S^1)^i$ is an $i$-dimensional compact torus.
D. Joyce (\cite{_Joyce_}) 
and independently Spindel et al. (\cite{_SSTvP_})
have shown that $G$ is equipped with a family of natural 
left-invariant hypercomplex structures. 
The corresponding complex structures were constructed
in 1950-ies by Bott and Samelson (\cite{_Samelson_}).

Consider the Killing metrics $(\cdot, \cdot)$
on $G$. Then $(\cdot, \cdot)$ is hypercomplex Hermitian,
with respect to the hypercomplex structure 
obtained by D. Joyce. Consider $G$ 
as a hypercomplex Hermitian manifold. 
Then $G$ is an HKT-manifold (\cite{_OP_}). 

\hfill

Denote the unit of $G$ by $e$. The
corresponting Lie algebra $\g$ is identified
with $T_eG$. Denote the structure
constants of $\g$ by  
\[ C^a_{b,c}:\; T_eG\times T_eG\arrow T_e G.
\]
Using the Killing form, we may identify
$T_e G$ and $T^*_eG$. Using this identification,
we consider $C^a_{b,c}$ as a 3-form
$\c T:\; T_eG\times T_eG\times T_eG\arrow \R$.
It is well known that $\c T$ is totally antisymmetric.
An exterior form on $T_eG$ can be uniquely extended
to a left invariant differential form on $G$.
Denote by $\goth T$ the differential form obtained from
$\c T$ this way. Then $\goth T$ is the torsion form of 
the HKT-connection $\nabla$ on $G$ (\ref{_HKT_Definition_}).
Moreover, $\nabla$ is the standard left-linvariant
torsion connection on the group $G$; a vector
field $v$ is parallel with respect to $\nabla$ 
if and only if $v$ is left-invariant.
 
Further details on this can be found in \cite{_Gra_Poon_}.

\hfill

The definition of an HKT-manifold is somewhat unwieldy,
as it uses the ``black box'' of the existence and
uniqueness of Bismut connection
(\ref{_Bismut_conne_exi_Theorem_}).
It is better to use the following theorem instead,
which is a version of a result from \cite{_Howe_Papado_}.

\hfill

\theorem\label{_HKT_via_diff_of_Omega_Theorem_}
(\cite{_Gra_Poon_}, Proposition 2)
Let $M$ be a quaternionic K\"ahler manifold. Given
a triple of induced complex
structures \[ I, J, K, \ \ I \circ J = - J \circ I = K,\] 
consider the corresponding non-degenerate (2,0)-form  
$\Omega\in \Lambda^{2,0}_I(M)$ \eqref{_sta_2,0_form_Equation_}.
Let $\6:\; \Lambda^{2,0}_I(M)\arrow \Lambda^{3,0}_I(M)$
be the Dolbeault differential, and $d$ the de Rham differential.
Then the following conditions are equivalent

\begin{description}
\item[(i)] $\6 \Omega =0$.
\item[(ii)] 
Let $\omega_I$, $\omega_J$, $\omega_K$ be the standard 
2-forms on $M$ \eqref{_sta_2,0_form_Equation_}.
The 3-forms $d\omega_I$, $d\omega_J$, $d\omega_K$ all
 have weight 1 with respect to the natural $SU(2)$-action
on $\Lambda^3(M)$.\footnote{This means that 
each of these forms belongs to a weight 1
$SU(2)$-subspace of $\Lambda^3(M)$.
An $SU(2)$-representation is of weight 1
if it is a direct sum of several copies
of a standard 2-dimensional irreducible
representation.}
\item[(iii)] $M$ is an HKT-manifold.

\end{description}

{\bf Proof:} 
The equivalence (i) $\Leftrightarrow$ (iii)
is proven in \cite{_Gra_Poon_}, Proposition 2. 
The implication (ii) $\Rightarrow$ (i) is clear,
because $\6 \Omega$ is a linear combination
of the 3-forms  $d\omega_J$, $d\omega_K$,
  $d^c\omega_J$, $d^c\omega_K$
(by $d^c$ we denote the ``twisted de Rham differential'',
$d^c := - I \circ d \circ I$). 
These four 3-forms all belong to the subrepresentation
of $\Lambda^3(M)$ generated by $d\omega_J$, $d\omega_K$,
hence by (ii) they have weight 1 with respect to $SU(2)$.
On the other hand, by \ref{_su(2)_action_explici_Lemma_}, 
a $(3,0)$-form $\6 \Omega\in \Lambda^{3,0}_I(M)$ has weight 3
if it is non-zero. Finally, an implication 
(iii) $\Rightarrow$ (ii) follows directly
from the Corollary 1 of \cite{_Gra_Poon_}
(see also \cite{_Howe_Papado_}).

One can easily deduce the proof of
\ref{_HKT_via_diff_of_Omega_Theorem_} 
using the quaternionic Dolbeault complex
(see \ref{_G_P_rheore_heuristic_Remark_}).
\endproof

\hfill

\remark
The standard $SU(2)$-action on $\Lambda^3(M)$
has weights 3 and 1. It is easy to see that the
weight 3 component is generated by
$\Lambda^{3,0}_I(M)$ for all induced
complex structures $I \in \H$. This explains the
implication
\ref{_HKT_via_diff_of_Omega_Theorem_} (i) $\Rightarrow$
\ref{_HKT_via_diff_of_Omega_Theorem_} (ii).

\hfill

\example
Let $M$ be a hypercomplex Hermitian manifold, of real dimension 4.
Since $\Lambda^3(M) \cong \Lambda^1(M)$, the group $SU(2)$
acts on $\Lambda^3(M)$ with weight 1. By 
\ref{_HKT_via_diff_of_Omega_Theorem_}, then,
$M$ is a HKT-manifold.


\section{The K\"ahler-de Rham superalgebra}
\label{_KdR_Section_}


In this Section,
we give a novel presentation of Kodaira relations,
which are interpreted as relations in a certain Lie superalgebra. 
These ideas
are commonplace in physics; for a physical
interpretations of the Lie superalgebra 
generated by de Rham and Dolbeault differentials
and the Hodge operators, see 
\cite{_FKS_}.

\subsection{Lie superalgebras}
\label{_LSA_Subsection_}

Let $A$ be a $\Z/2\Z$-graded vector space,
\[ A = A^\ev \oplus A^\odd.\] We say that $a\in A$
is {\bf pure} if $a$ belongs to $A^\ev$ or $A^\odd$.
For a pure element $a\in A$, we
write $\tilde a =0$ if $a\in A^\ev$,
and $\tilde a =1$ if $a\in A^\odd$.
Consider a bilinear operator
\[ [\cdot, \cdot]:\; A \times A \arrow A,\]
called {\bf supercommutator}. Assume that 
$[\cdot, \cdot]$ is {\bf graded an\-ti-\-com\-mu\-ta\-ti\-ve},
that is, satisfies
\[ [ a, b] = - (-1)^{\tilde a \tilde b}  [ b,a ]
\]
for pure $a, b \in A$. Assume, moreover, that
$[\cdot, \cdot]$ is compatible with the grading:
the commutator $[ a, b]$ is even when both $a$, $b$ are even 
or odd, and odd if one of these elements is odd and another 
is even. We say that $A, [\cdot, \cdot]$ is a Lie superalgebra,
if the following identity (called {\bf the graded Jacoby
identity}) holds, for all pure elements $a, b, c\in A$
\begin{equation}
[ a, [ b, c]] = [ [ a, b], c] + (-1)^{\tilde a \tilde b}[b, [ a,c ]].
\end{equation}
Up to a sign, this is a usual Jacoby identity.

Every reasonable property of Lie algebras has a 
natural analogue for Lie superalgebras; the rule of thumb,
every time one would exchange two elements $a$ and $b$,
one adds a multiplier $(-1)^{\tilde a \tilde b}$.

In classical terms, one interprets the Lie superalgebras as follows:
$A^\ev$ is a usual Lie algebra, $A^\odd$ is an $A^\ev$-module,
equipped with a bilinear symmetric $A^\ev$-invariant pairing
\[ [\cdot, \cdot]:\; A^\odd \times A^\odd \arrow A^\ev\]
which satisfies the (super-)Jacoby identity 
\[ 
[ a, [ b, c]] = [ [ a, b], c] + [ [ a,c ], b]
\]
for all $a, b, c \in A^\odd$.

\hfill

\example 
Let $V= V^\ev\oplus V^\odd$ be a $\Z/2\Z$-graded vector space,
and $\End V$ its endomorphism space, equipped with induced grading.
We define a supercommutator in $\End V$ by the formula
\[ [a, b] = ab - (-1)^{\tilde a \tilde b}ba
\]
It is easy to check that $\left(\End V, [\cdot, \cdot]\right)$
is a Lie superalgebra. 

\subsection{Kodaira relations and the K\"ahler-de Rham superalgebra}

Let $M$ be a K\"ahler manifold. Consider $\Lambda^*(M)$ as a graded
vector space. The differentials $d, d^c:= - I \circ d \circ I$
can be interpreted as odd elements in $\End(\Lambda^*(M))$,
and the Hodge operators $L, \Lambda, H$ as even elements. 
To avoid confusion, we denote the supercommutator of odd
elements as $\{\cdot,\cdot\}$. In terms of the associative
algebra, $\{a,b\} = ab +ba$.
Let $d^*:=  [\Lambda, d^c]$, $(d^c)^*:= -[\Lambda, d]$.
The usual Kodaira relations can be stated as follows

\begin{equation}\label{_Kodaira_Equations_}
\begin{aligned} 
 {} & [ L, d^*] = - d^c,\ \ \ \  [ L, (d^c)^*] = d, \ \ \ \ 
\{ d, (d^c)^*\} = \{ d^*, d^c\} =0, \\
 {} & \{ d, d^c\} =  \{ d^*, (d^c)^*\} =0,\ \ \ \  
\{ d, d^*\} = \{ d^c, (d^c)^*\} = \Delta,
\end{aligned}
\end{equation}
where $\Delta$ is the Laplace operator,
 commuting with $L, \Lambda, H$, and
$d$, $d^c$.

\hfill

\definition \label{_KdR_Definition_}
Let $M$ be a K\"ahler manifold.
Consider the Lie superalgebra $\g \subset \End(\Lambda^*(M))$
generated by the even elements $L, \Lambda, H$, and
the odd elements $d$, $d^c$ Then $\g$ is called
{\bf the K\"ahler-de Rham superalgebra} associated with $M$.

\hfill

Using \eqref{_Kodaira_Equations_}, 
is easy to see that
$\g$ is in fact independent form $M$
(see Subsection \ref{_KdR_axioma_Subsection_}).

\hfill

This Lie superalgebra was studied from the physicists'
point of view in \cite{_FKS_}.

\subsection{Axiomatics of the K\"ahler-de Rham superalgebra}
\label{_KdR_axioma_Subsection_}

The K\"ahler-de Rham superalgebra has the following
formal interpretation. 

\hfill

\theorem\label{_Kah_deRha_formally_Theorem_}
Let $\g$ be a Lie superalgebra, and $\goth a\subset \g^\ev$ its
even subalgebra isomorphic to 
$\goth {sl}(2)$,  and generated by the
standard $\goth {sl}(2)$-triple 
$\langle L, \Lambda, H\rangle\in \g^\ev$.
Consider odd vectors $d, d^c \subset \g^{odd}$.
Assume that $[ H, d] = d$, $[H, d^c] = d^c$,
and $[ L, d] = [L, d^c] =0$, so that $d$ and $d^c$ generate
weight-1-representations of $\goth a \cong \goth {sl}(2)$.
Assume, moreover, that $\g$ is generated as
a superalgebra by $(L, \Lambda, H, d, d^c)$
Finally, assume that
\begin{equation} \label{_d,d_d^c_commute_Equation_}
  \{ d, d\} = \{d^c, d^c\} = \{ d, d^c\} =0, 
\end{equation}
where $\{\cdot, \cdot\}$ denotes the 
super-commutator of odd elements.
Consider the vectors
\[ d^* := [ \Lambda, d^c] ,\ \ (d^c)^* := -[ \Lambda, d].
\]
Then
\begin{description}
\item[(i)] We have
\begin{equation}\label{_Laplace_in_KdR_formal_Equation_} 
  \{ d, d^*\} = \{d^c, (d^c)^*\}.
\end{equation}
Denote this supercommutator by
\[ \Delta:= \{ d, d^*\} = \{d^c, (d^c)^*\}.
\]

\item[(ii)] The Lie superalgebra $\g$ is 8-dimensional,
and spanned by \[ L, \Lambda, H, d, d^c, d^*, (d^c)^*, \Delta.\]
Moreover, $\g$ is naturally isomorphic to the 
K\"ahler-de Rham superalgebra (\ref{_KdR_Definition_}).
\end{description}

{\bf Proof:} 
The action of the $\goth{sl}(2)$-triple
$\langle L, \Lambda, H\rangle$
on \[ V:= \langle d, d^c, d^*, (d^c)^*\rangle\subset \g\]
is transparent, because $V$ is a direct sum of two 
weight 1 representations, $\langle d, (d^c)^*\rangle$
and $\langle d^c, d^*\rangle$. 
To prove \ref{_Kah_deRha_formally_Theorem_},
we have to show that the commutator relations in the
4-dimensional odd space $\langle d, d^c, d^*, (d^c)^*\rangle\subset \g$
are consistent with \eqref{_Kodaira_Equations_},
and to show that $\Delta$ belongs to the center of $\g$.

We start with proving \eqref{_Laplace_in_KdR_formal_Equation_}.
Consider the relation
\begin{equation}\label{_d_d^c_commute_Equation_}
\{ d, d^c\} =0
\end{equation}
(see \eqref{_d,d_d^c_commute_Equation_}).
Since the super-commutator in $V$ is $\goth a$-invariant,
\eqref{_d_d^c_commute_Equation_} implies
\[ [ \Lambda, \{ d, d^c\}] = 
  \{ [ \Lambda,d], d^c\} +\{ d, [ \Lambda,d^c]\} =0
\]
which means
\[ \{ (d^c)^*, d^c\} - \{ d, d^*\} =0. \]
This proves \ref{_Kah_deRha_formally_Theorem_} (i). 
Remaining commutation relations between
$d, d^c, d^*, (d^c)^*$ are obtained in a similar way,
as follows.

Acting on $\{d, d\}=0$ with $[\Lambda, \cdot]$,
we obtain
\[ \{ [\Lambda, d], d\} =0, \]
that is, $\{ (d^c)^*, d\}=0$. Similarly, 
$\{d^c, d^c\}=0$ acted on by $[\Lambda, \cdot]$
implies $\{d^*, d^c\}=0$. Since the vectors
$d, d^c$ have weight 1 with respect to the
$\goth{sl}(2)$-action $L, \Lambda, H$, 
we have
\[ [\Lambda, [ \Lambda, d]] = [\Lambda, [ \Lambda, d^c]] =0. \]
In other words, 
\begin{equation}\label{_Lambda_on_d^*_Equation_} 
  [\Lambda, d^*] = [\Lambda, (d^c)^*] =0.
\end{equation}
Acting with $[\Lambda, \cdot]$ on both sides of 
$\{d^*, d^c\}=0$, and using
\eqref{_Lambda_on_d^*_Equation_}, we obtain
\[ 0= [\Lambda, \{d^*, d^c\}] = \{[\Lambda, d^*], d^c\} +
   \{d^*, [\Lambda,d^c]\} = \{d^*, [\Lambda,d^c]\} = \{d^*,d^*\}
\]
Using $\{ (d^c)^*, d\}=0$ instead of $\{d^*, d^c\}=0$
and applying the same argument, obtain $\{ (d^c)^*, (d^c)^*\}=0$.
Finally, if we apply $[\Lambda, \cdot]$  to
\[ \{ (d^c)^*, d^c\} - \{ d, d^*\} =0, \]
we shall have
\[ \{ (d^c)^*, d^*\} + \{ (d^c)^*, d^*\} =0,\]
that is, 
\begin{equation}\label{_d^c^*d^*_commute_Equation_}
\{ (d^c)^*, d^*\} =0.
\end{equation}
We have found that pairwise super-commutators
of the vectors $d, d^c, d^*, (d^c)^*$
are equal zero, except
\[ \{ d, d^*\} = \{d^c, (d^c)^*\} = \Delta.
\]
We obtained the Kodaira identities 
\eqref{_Kodaira_Equations_}.
To prove \ref{_Kah_deRha_formally_Theorem_},
it remains to show that the ``Laplacian'' 
$\Delta$ commutes with 
\[ d, d^c,  L, \Lambda, H.\]

Writing $\Delta$ as $\{ d, d^*\}$ and using the commutation
relations between $d^c$, $(d^c)^*$  and $d, d^*$, we find that
$\Delta$ commutes with $d^c$  and $(d^c)^*$. Similarly,
using the formula $\Delta= \{d^c, (d^c)^*\} $ and
the fact that $d$, $d^*$ commutes with $d^c$, $(d^c)^*$,
we find that $d$, $d^*$ commutes with $\Delta$. The operator
$H$ commutes with $\Delta$ because $H$ multiplies
$d$ by 1 and $d^*$ by $-1$. Finally, $L$ and $\Lambda$
commute with $\Delta$ because $L$ maps $d$ to zero
and $d^*$ to $d^c$, hence
\[ [ L, \{ d, d^*\}] = \{ [ L, d], d^*\} + \{ d,  [ L,d^*]\}=
   \{ d, d^c\};
\]
this commutator is equal zero by 
\eqref{_d_d^c_commute_Equation_}. Similarly,
\[
[ \Lambda, \{ d, d^*\}] = \{ [ \Lambda, d], d^*\} + \{ d,  [ \Lambda,d^*]\}=
   -\{ d^*, (d^c)^*\};
\]
this is equal zero by \eqref{_d^c^*d^*_commute_Equation_}.
We have proven \ref{_Kah_deRha_formally_Theorem_}.
\endproof

\subsection{K\"ahler-de Rham superalgebra and the odd Heisenberg superalgebra}

\definition
Let $\g$ be the Lie superalgebra defined
by generators and relations as in 
\ref{_Kah_deRha_formally_Theorem_}. Then $\g$ 
is called {\bf the K\"ah\-ler-\-de Rham superalgebra}.

\hfill

The K\"ahler-de Rham superalgebra can be described
explicitly as follows. Let $V$ be a weight 1 representation of
$\goth{sl}(2)$. Consider the space $V\oplus V^*$ with the
natural $\goth{sl}(2)$-invariant 
bilinear symmetric form $\langle \cdot, \cdot, \rangle$.
Consider the graded space 
$\goth{h}:= V\oplus V^*\oplus \C c$,
with $\C c$ 1-dimensional even space and
$V\oplus V^*$ odd. We define a supercommutator
on $\goth h$ in such a way that
\[ \{v_1, v_2\}:= \langle v_1, v_2 \rangle c,
\ \ \text{for all}\ \ v_1, v_2 \in V\oplus V^*
\]
and $c$ is a central element.
The Lie superalgebra $\goth h$ is called
{\bf odd Heisenberg algebra associated with 
$V\oplus V^*$}. By construction,
$\goth{sl}(2)$ acts in $\goth h$ by 
Lie superalgebra automorphisms. 

\hfill

\proposition\label{_Kah_deRham_explici_Proposition_}
In the above notation, consider the semidirect product
$\goth{sl}(2)\rightthreetimes \goth{h}$. Then
$\goth{sl}(2)\rightthreetimes \goth{h}$ is
isomorphic to the K\"ahler-de Rham Lie superalgebra
$\g$.

\hfill

{\bf Proof:} Consider the odd 
4-dimensional subspace 
$W\subset \g$ spanned by
$d$, $d^c$, $d^*$, $(d^c)^*$.
By Kodaira idenitities, the
 commutator of vectors $v, v'\in W$
is proportional to the Laplacian.
Let 
\[ \langle v, v' \rangle:= \frac{\{ v, v'\}}{\Delta}.
\]
Clearly, $W$ generates a 5-dimensional subalgebra
${\goth k}\subset \g$ which is naturally isomorphic to $\goth h$.
The Hodge operators $L, \Lambda, H$
commute with the Laplacian, and therefore
preserve the bilinear form 
$\langle v, v' \rangle$. Therefore,
the $\goth{sl}(2)$-algebra generated by
$L, \Lambda, H$ acts on
${\goth k}$ in the same way as
$\goth{sl}(2)$ acts on $\goth{h}$.
Finally, $\g$ is by construction
a semidirect product of $\goth{sl}(2)$ and
${\goth k}$. \endproof


\section{Quaternionic Dolbeault complex}
\label{_q.D._Section_}


The quaternionic cohomology is a well 
known subject, introduced by M. Capria and S. Salamon
(\cite{_Capria-Salamon_}). 
Here we give an
exposition of quaternionic cohomology and quaternionic
Dolbeault complex for hypercomplex manifolds.
We follow \cite{_V:projective_}.

\subsection{Quaternionic Dolbeault complex: the definition}

Let $M$ be a hypercomplex manifold, and 

\[ \Lambda^0 M \stackrel d\arrow \Lambda^1 M 
   \stackrel d\arrow \Lambda^2 M \stackrel d\arrow ...
\]
its de Rham complex. Consider the natural action of $SU(2)$ on 
$\Lambda^*M$.  Clearly, $SU(2)$ acts on 
$\Lambda^iM$, $i\leq \frac 1 2 \dim_\R M$
with weights $i, i-2, i-4, \dots$

We denote by $\Lambda^i_+$ the maximal $SU(2)$-subspace 
of $\Lambda^i$, on which $SU(2)$ acts with weight $i$. 

\hfill

The following linear-algebraic lemma allows one to compute
$\Lambda^i_+$ explicitly

\hfill

\lemma\label{_Lambda_+_explicit_Lemma_}
In the above assumptions, let $I$ be an induced complex structure,
and $\H_I$ the quaternion space, considered as a 2-dimensional
complex vector space with the complex structure induced by $I$.
Denote by $\Lambda^{p,0}_I(M)$ the space of $(p,0)$-forms on $(M,I)$. 
The space $\H_I$ is equipped with a
natural action of $SU(2)$. Consider $\Lambda^{p,0}_I(M)$
as a representation of $SU(2)$, with trivial group action. 
Then, there is a canonical isomorphism
\begin{equation}\label{_Lambda_+_explicit_Equation_}
\Lambda^p_+(M) \cong S^p_\C \H_I \otimes_\C \Lambda^{p,0}_I(M),
\end{equation}
where $S^p_\C \H_I$ denotes a $p$-th symmetric power of $\H_I$.
Moreover, the $SU(2)$-action on $\Lambda^p_+(M)$ is compatible with
the isomorphism \eqref{_Lambda_+_explicit_Equation_}.

\hfill

{\bf Proof:} This is \cite{_V:projective_}, Lemma 8.1.
\endproof

\hfill

Consider an $SU(2)$-invariant decomposition
\begin{equation}\label{_decompo_Lambda_to_Lambda_+_Equation_}
\Lambda^p(M) = \Lambda^p_+(M)\oplus V^p,
\end{equation}
where $V^p$ is the sum of all $SU(2)$-subspaces 
of $\Lambda^p(M)$ of weight less than $p$. Using the decomposition
\eqref{_decompo_Lambda_to_Lambda_+_Equation_},
we define the quaternionic Dolbeault differential 
$d_+:\; \Lambda^*_+(M)\arrow \Lambda^*_+(M)$ 
as a composition of de Rham differential and projection of to
$\Lambda^*_+(M) \subset \Lambda^*(M)$. Since the 
de Rham differential cannot increase the $SU(2)$-weight
of a form more than by 1, $d$ preserves
the subspace $V^*\subset \Lambda^*(M)$.
Therefore, $d_+$ is a differential in 
$\Lambda^*_+(M)$.

\hfill

\definition
Let
\[ \Lambda^0 M \stackrel {d_+}\arrow \Lambda^1 M 
   \stackrel {d_+}\arrow \Lambda^2_+ M \stackrel {d_+}\arrow 
   \Lambda^3_+ M \stackrel {d_+}\arrow...
\]
be the differential graded algebra constructed above\footnote{We 
identify $\Lambda^0 M$ and $\Lambda^0_+ M$,
$\Lambda^1 M$ and $\Lambda^1_+ M$.}.
It is called {\bf the quaternionic Dolbeault complex},
or {\bf Salamon complex}.

\subsection[Hodge decomposition for the quaternionic Dolbeault complex]{Hodge decomposition for the quaternionic  \\ Dolbeault complex}
\label{_Hodge_Dolbeault_Subsection_}

Let $M$ be a hypercomplex manifold, and $I$ an induced complex
structure.  Consider the operator 
$adI:\; \Lambda^*(M) \arrow \Lambda^*(M)$
mapping a $(p,q)$-form $\eta$ to $\1(p-q)\eta$.
By definition, $ad I$ belongs to the Lie algebra $\goth{su}(2)$
acting on $\Lambda^*(M)$ in the standard way. 
Therefore, $ad I$ preserves the subspace
$\Lambda^*_+(M) \subset \Lambda^*(M)$. 
We obtain the Hodge decomposition
\[ \Lambda^*_+(M) = \oplus_{p,q} \Lambda^{p,q}_{+,I}(M). \]

\hfill

\definition\label{_Hodge_for_q.D._Definition_}
The decomposition
\[ \Lambda^*_+(M) = \oplus_{p,q} \Lambda^{p,q}_{+,I}(M)
\]
is called {\bf the Hodge decomposition for the 
quaternionic Dolbeault complex}.

\hfill

The following claim is trivial

\hfill

\claim \label{_Lambda_+^{p,0}_Claim_}
Given a hypercomplex manifold $M$ and 
an induced complex structure $I$, the following
subspaces of $\Lambda^*(M)$ coinside:
\[ 
  \Lambda^{p,0}_{+,I}(M) = \Lambda^{p,0}_{I}(M),
\]
where $\Lambda^{p,0}_{I}(M)$ denotes the space of all $(p,0)$-forms.

\hfill

{\bf Proof:} Immediately follows from
\ref{_Lambda_+_explicit_Lemma_}. 
\endproof

\subsection{Quaternionic Dolbeault bicomplex: explicit description}
\label{_qD-bico_Subsection_}

Let $M$ be a hypercomplex manifold, $I$ an induced comlex structure,
and $I, J, K\in \H$ the standard triple of induced complex structures. 
Clearly, $J$ acts on the complexified co tangent space
$\Lambda^1M\otimes \C$ mapping $\Lambda_I^{0,1}(M)$
to $\Lambda_I^{1,0}(M)$. 
Consider a differential operator 
\[ \6_J:\; C^\infty(M)\arrow \Lambda_I^{1,0}(M),\]
mapping $f$ to $J(\bar\6 f)$, where 
$\bar\6:\; C^\infty(M)\arrow \Lambda_I^{0,1}(M)$
is the standard Dolbeault differential on a K\"ahler
manifold $(M,I)$. We extend $\6_J$
to a differential
\[ 
   \6_J:\; \Lambda_I^{p,0}(M)\arrow \Lambda_I^{p+1,0}(M),
\]
using the Leibniz rule.

\hfill

\proposition \label{_d_+_Hodge_6_J_Proposition_}
Let $M$ be a hypercomplex manifold, $I$ an induced 
complex structure, $I, J, K$ the standard 
basis in quaternion algebra, and
\[ \Lambda^*_+(M)= \oplus_{p,q} \Lambda^{p,q}_{I,+}(M) \]
the Hodge decomposition of the quaternionic Dolbeault
complex. Then there exists a canonical isomorphism
\begin{equation}\label{_Lambda_+_Hodge_deco_expli_Equation_} 
   \Lambda^{p,q}_{I, +}(M)\cong \Lambda^{p+q, 0}_I(M). 
\end{equation}
Under this identification, the quaternionic Dolbeault differential
\[ d_+:\; \Lambda^{p,q}_{I, +}(M)\arrow 
   \Lambda^{p+1,q}_{I, +}(M)\oplus \Lambda^{p,q+1}_{I, +}(M)
\]
corresponds to a sum
\[ \6 \oplus \6_J:\; \Lambda^{p+q, 0}_{I}(M) \arrow
   \Lambda^{p+q+1, 0}_{I}(M)\oplus \Lambda^{p+q+1,0}_{I}(M).
\]

{\bf Proof:} This is Proposition 8.13 of \cite{_V:projective_}.
\endproof

\hfill

The statement of 
\ref{_d_+_Hodge_6_J_Proposition_}
can be represented 
by the following diagram

\begin{equation}\label{_bicomple_XY_Equation}
\begin{minipage}[m]{0.85\linewidth}
{\tiny $
\xymatrix @C+1mm @R+10mm@!0  { 
  {} &{} & \Lambda^0_+(M) \ar[dl]^{d'_+} \ar[dr]^{d''_+} 
   {} &{} & {} &{} & {} &{} & \Lambda^{0,0}_I(M) \ar[dl]^{ \6} \ar[dr]^{ \6_J}
   {} &{} &  \\
 {} & \Lambda^{1,0}_+(M) \ar[dl]^{d'_+} \ar[dr]^{d''_+} {} &
 {} & \Lambda^{0,1}_+(M) \ar[dl]^{d'_+} \ar[dr]^{d''_+}{} &{} & 
\text{\large $\cong$} {} &
 {} &\Lambda^{1,0}_I(M)\ar[dl]^{ \6} \ar[dr]^{ \6_J}{} &  {} &
 \Lambda^{1,0}_I(M)\ar[dl]^{ \6} \ar[dr]^{ \6_J}{} &\\
 \Lambda^{2,0}_+(M) {} &{} & \Lambda^{1,1}_+(M) 
   {} &{} & \Lambda^{0,2}_+(M){} & \ \ \ \ \ \ {} & \Lambda^{2,0}_I(M){} & {} & 
\Lambda^{2,0}_I(M) {} & {} & \Lambda^{2,0}_I(M) \\
}
$
}
\end{minipage}
\end{equation}
where $d_+= d'_+ + d''_+$ is the Hodge decomposition of 
the quaternionic Dolbeault differential.

\hfill

\remark 
Since $d_+^2=0$, we obtain that 
$\6$ and $\6_J$ anti-commute:
\begin{equation}\label{_6_6_J_commute_Equation_}
\{\6, \6_J\} =0.
\end{equation}

\hfill

Further on in this paper, only
\eqref{_6_6_J_commute_Equation_} will be used;
however, the analogy between the 
quaternionic Dolbeault complex and the de Rham algebra
of a K\"ahler manifold is implicit
in our construction. 

\hfill

The definition of an HKT-manifold can be reformulated
in terms of a quaternionic Dolbeault complex as follows.

\hfill

\theorem \label{_HKT_in_terms_of_d_+_Theorem_}
Let $M$ be a hypercomplex Hermitian manifold,
$I$ an induced complex structure, $\omega_I$
the corresponding 2-form, and $d_+$ the quaternionic
Dolbeault differential. Then $\6 \Omega=0$
if and only if $d_+\omega_I =0$. 

\hfill

{\bf Proof:} 
Let $I,J,K$ be the standard basis in quaternions.
Consider $\omega_I$ as an element
in $\Lambda^{1,1}_{I,+}(M)$. Using 
\ref{_d_+_Hodge_6_J_Proposition_},
we identify $\Lambda^{1,1}_{I,+}(M)$
and $\Lambda^{2,0}_{I,+}(M)$. Under
this identification, $\omega_I$
corresponds to the standard $(2,0)$-form $\Omega$
(see \cite{_V:projective_}).
Therefore, $\omega_I$ is $d_+$-closed
if and only if $\Omega$ is $d_+$-closed.

Clearly,
\[ \6_J \Omega = J(\overline{\6\Omega}) \]
(see, e.g. \ref{_commu_L_Omega_6_J_Proposition_}),
hence $\6 \Omega=0$ if and only if $\6_J \Omega=0$.
By \ref{_d_+_Hodge_6_J_Proposition_}, we find that
$\6 \Omega=0$ if and only if
$d_+\Omega=0$. This proves
\ref{_HKT_in_terms_of_d_+_Theorem_}.
\endproof

\hfill

\remark
By \ref{_HKT_via_diff_of_Omega_Theorem_},
$M$ is an HKT-manifold if and only if 
$\6 \Omega=0$. This allows one to interpret
\ref{_HKT_in_terms_of_d_+_Theorem_}
in terms of the HKT-condition: a manifold
$M$ is HKT if and only if $d_+\omega_I =0$.
This statement is analogous to the 
definition of a K\"ahler manifold
in terms of a symplectic form:
 ``a complex Hermitian manifold
is K\"ahler if and only if the standard
$(1,1)$-form $\omega$ is closed''.

\hfill

\remark \label{_G_P_rheore_heuristic_Remark_}
One can easily deduce 
\ref{_HKT_via_diff_of_Omega_Theorem_} 
(the HKT conditions)
from \ref{_HKT_in_terms_of_d_+_Theorem_}.
Consider a hypercomplex Hermitian manifold
which satisfies $\6 \Omega=0$. Then
$d_+\Omega=0$ and $d_+\omega_I=0$.
Therefore, for $\omega = \omega_I, \omega_J, \omega_K$,
the differential $d\omega$ is of weight 1 with respect
to the $SU(2)$-action.  
This is the condition
(ii) of \ref{_HKT_via_diff_of_Omega_Theorem_}.
In other words, 
\ref{_HKT_in_terms_of_d_+_Theorem_}
allows us to check that \ref{_HKT_via_diff_of_Omega_Theorem_} 
(i) $\Leftrightarrow$ \ref{_HKT_via_diff_of_Omega_Theorem_} 
(ii). 

The implication
\ref{_HKT_via_diff_of_Omega_Theorem_} 
(ii) $\Rightarrow$ \ref{_HKT_via_diff_of_Omega_Theorem_} 
(iii) is standard. The torsion forms of the connections
$\nabla_I$, $\nabla_J$, $\nabla_K$ are equal to
$Id \omega_I$,  $Jd \omega_J$, $Kd \omega_K$.
Suppose that
\begin{equation}\label{_torsion_forms_Equation_}
Id \omega_I=Jd \omega_J=Kd \omega_K,
\end{equation}
that is, the torsion forms of $\nabla_I$, $\nabla_J$, $\nabla_K$
are equal.
Substracting  the torsion form from $\nabla_I$,
$\nabla_J$, $\nabla_K$, we obtain three
torsion-free orthogonal connections on $M$.
These connections are equal to Levi-Civita 
connection, because Levi-Civita 
connection is the unique orthogonal torsion-free 
connection. This implies $\nabla_I=\nabla_J=\nabla_K$.
Therefore, a hypercomplex Hermitian manifold is HKT
if and only if \eqref{_torsion_forms_Equation_}
holds. When $Id \omega_I$,  $Jd \omega_J$, $Kd \omega_K$
all have weight one with respect to $SU(2)$,
\eqref{_torsion_forms_Equation_} follows from
elementary calculation; converse is also clear.


\section{The hyperk\"ahler-de Rham superalgebra}
\label{_hk-DR_Section_}


\subsection{The $\goth{so}(1,4)$-action on $\Lambda^*(M)$}

Let $M$ be a hyperk\"ahler manifold, $I$, $J$, $K$
induced complex structures giving a standard basis
in quaternions, and $\omega_I$, $\omega_J$,
$\omega_K$ the corresponding 2-forms (\ref{_induced_c_str_Definition_}).
Denote by $L_I, L_J, L_K:\; \Lambda^*(M) \arrow \Lambda^{*+2}(M)$ 
the operators of exterior multiplication by $\omega_I$, $\omega_J$,
$\omega_K$. Let 
\[ 
\Lambda_I:= * L_I*, \Lambda_J:= * L_J*, \Lambda_K:= * L_K*
\]
be the Hermitian adjoint operators, and $H$ the 
standard Hodge operator; we have three $\goth{sl}(2)$-triples
$\langle L_I, \Lambda_I, H\rangle$, 
$\langle L_J, \Lambda_J, H\rangle$, 
$\langle L_K, \Lambda_K, H\rangle$
acting on $\Lambda^*(M)$. 

\hfill

Given an induced complex structure $L\in \H$, consider
an operator 
\[ adL:\; \Lambda^{p,q}_L(M) \arrow \Lambda^{p,q}_L(M),
   \ \ \ \eta \arrow \1 (p-q) \eta. 
\]
Clearly, $adL$ is the Lie algebra element
corresponding to the $U(1)$-action which induces
the Hodge decomposition. Therefore, the endomorphism
$ad L$ belongs to the Lie algebra of the
standard $SU(2)$-action on $H^*(M)$; moreover,
this Lie algebra is generated by $ad L$,
for $L = I, J, K$.

\hfill

The following theorem is implied by an easy linear-algebrain
computation. 

\hfill

\theorem\label{_SO(4,1)_Theorem_}
Consider the 10-dimensional vector space space 
$\goth a \subset \End(\Lambda^*(M))$ generated by
$H, L_X, \Lambda_X, ad X$, for $X= I, J, K$.
Then $\goth a$ is a Lie subalgebra of
$\End(\Lambda^*(M))$. Moreover, $\goth a$
is isomorphic to $\goth{so}(1,4)$.

\hfill

{\bf Proof:} See \cite{_Verbitsky:SO(5)_}. \endproof

\subsection{Hyperk\"ahler-de Rham superalgebra: the definition}

Consider the 2-dimensional
quaternionic vector space $W$ equipped with a quaternionic
Hermitian metrics of signature $+,-$. 
Denote the group of quaternionic-linear
isometries of $W$ by $SU(1,1, \H)$.
The group $SU(1,1, \H)$ acts on
$W$ in a tautological way.\
Using the classical isomorphism
\[ SU(1,1, \H)\cong Spin(1,4), \]
we may consider $W$ as a representation of
$\goth a =\goth{so}(1,4)$.

Assume that the grading of $W$ is odd.
Let $\c H\; :=  W \oplus \R D$ be a direct
sum of $W$ and an even 1-dimensional vector space.
Consider the map 
\[ 
   (\cdot, \cdot)D:\; W\times W \arrow  \R D
\] 
putting $v, v' \in W$
to $(v, v')_D\subset \c H$, where $(\cdot, \cdot)$
is the Hermitian form on $W$, and $D$ the
generator of $\R D\subset \c H$.

We introduce a structure of a Lie superalgebra on $\c H$
as follows: for odd elements $v, v' \in W$, we have  
\[ \{v, v'\} = (v, v')_D.
\]
and a supercommutator of $D$ with anything is zero.
Such an algebra is called {\bf the odd Heisenberg algebra};
it is a graded version of a usual Heisenberg algebra.

Clearly, $\goth a =\goth{so}(1,4)$ acts on 
$\c H$ by automorphisms.

\hfill

\definition\label{_hk-dR_Definition_}
Let $\goth a \rtimes \c H$ be a semidirect product of
$\goth a$ and $\c H$. Then $\goth a \rtimes \c H$ 
{\bf a hyperk\"ahler-de Rham superalgebra}.

\hfill

\theorem \label{_hk-dR-acts_on_forms_Theorem_}
(\cite{_FKS_})
Let $M$ be a hyperk\"ahler manifold, $\goth a \cong\goth{so}(1,4)$
the Lie algebra constructed in \ref{_SO(4,1)_Theorem_}
and acting on the differential forms,
$\goth a\subset \End(\Lambda^*(M))$.
Consider the Lie superalgebra
$\g \subset \End(\Lambda^*(M))$,
generated by the even subspace
$\goth a\subset \End(\Lambda^*(M))$
and the odd vector $d\in \End(\Lambda^*(M))$,
where $d:\; \Lambda^*(M) \arrow \Lambda^{*+1}(M)$ 
is the de Rham differential. Then $\g$ is naturally
isomorphic to the hyperk\"ahler-de Rham superalgebra
(\ref{_hk-dR_Definition_}).

\hfill

{\bf Proof} Let ${I_1}\in \H$ be an induced complex structure.
Consider the differential
$d_{I_1}:= - {I_1} \circ d \circ {I_1}$.\footnote{In the K\"ahler
situation, this differential was denoted by $d^c$.
Over a hyperk\"ahler manifold, one needs to emphasize
the dependence from the choice of the induced
complex structure ${I_1}$.}. Since $d_{I_1} = \1 [ ad{I_1}, d]$,
the differential $d_{I_1}$ belongs to $\g$. 
Let $d^*_{I_1}:= - * d_{I_1} *$. 
By Kodaira identities
(see \eqref{_Kodaira_Equations_}), we have
\[ d^*_{I_1} = -[ \Lambda_{I_1}, d]. \]
Therefore, the vectors $d^*_{I_1}$ belong to 
$\g$.

Let $W\subset \g$ be the 8-dimensional odd subspace
spanned by \[ d, d_I, d_J, d_K, d^*, d^*_I, d^*_J, d^*_K. \]
Using the Kodaira relations, it is easy to check that
all pairwise commutators of the 
vectors $\langle d, d_I, d_J, d_K, d^*, d^*_I, d^*_J, d^*_K\rangle$
vanish, except $\{d, d^*\}$, $\{d_I, d_I^*\}$, 
$\{d_J, d_J^*\}$, $\{d_K, d_K^*\}$. Moreover, we have
\begin{equation}\label{_Laplace_hk_Equation_}
\{d, d^*\}=\{d_I, d_I^*\}=\{d_J, d_J^*\}= \{d_K, d_K^*\}=\Delta, 
\end{equation}
where $\Delta = d d^* + d^* d$ is the usual 
Laplace operator on the differential forms. 

By the Kodaira relations,
 the Laplacian $\Delta$ commutes with $d$, $d_{I_1}$, 
$L_{I_1}$ and $\Lambda_{I_1}$, for any induced K\"ahler 
structure $I_1$. Therefore, $\Delta$ belongs to the centre of $\g$.

For any two vectors $v, v' \in W$, the super-commutator
$\{ v, v'\}$ is proportional to the Laplacian.
This allows one to speak of the quotient
$\frac{\{ v, v'\}}{\Delta}\in \R$.
Consider the pairing
\[ (\cdot, \cdot):\; W\times W \arrow \R \]
mapping $v, v'$ to $\frac{\{ v, v'\}}{\Delta}$.
Since $\goth a \cong \goth{so}(1,4)$
preserves the Laplacian, $\goth a$
acts on $W$ preserving the pairing
$(\cdot, \cdot)$. An elementary calculation is used to check that
$(\cdot, \cdot)$ has signature $(+,+,+,+,-,-,-,-)$.
The Lie algebra  $\goth a =\goth{su}(1,1, \H)$
acts on $W$ as on its fundamental representation;
this identifies $W$ and $\H \oplus \H$.
Since $\Delta$ lies in the center of $\g$, the 
subalgebra $\c H:= W \oplus \R \Delta \subset \g$
is isomorphic to the odd Heisenberg algebra;
this sublalgebra is preserved by 
$\goth a =\goth{so}(1,4)\subset \g$,
and $\g$ is generated by $\goth a$
and $\c H$. Therefore, $\c H$ is an ideal of $\g$, and
$\g $ is a semidirect product of $\c H$ and $\goth a$.

We have proven \ref{_hk-dR-acts_on_forms_Theorem_}. \endproof

\subsection{Dolbeault complex and 
the hyperk\"ahler-de Rham superalgebra}
\label{_Dolbeault_hype_KdR_Subsection_}


Let $M$ be a hyperk\"ahler manifold, $I, J, K$
the standard triple of induced complex structures,
$I\circ J = - J \circ I = K$, and 
$\g = \goth a \rtimes \c H$ the
K\"ahler-de Rham superalgebra acting on $\Lambda^*(M)$.

\hfill

Let $\Omega= \frac 1 2 (\omega_J + \1\omega_K)$
be the standard $(2,0)$-form on $(M,I)$,
and $L_\Omega$ the operator
of multiplication by $\Omega$. Consider
the Hermitian adjoint operator 
$\Lambda_\Omega = *L_\Omega*$,
and let $H_\Omega$ be their commutator. 
Using \ref{_SO(4,1)_Theorem_}, it is easy to
check that $\langle L_\Omega,\Lambda_\Omega, H_\Omega\rangle$
is an $\goth{sl}(2)$-triple (see \cite{_Verbitsky:SO(5)_}).
Moreover, $H_\Omega$ maps $\eta\in \Lambda^{p, q}_I(M)$
to $(n-p)\eta$, where $n = \dim_\H M$.

\hfill

The main result of this Subsection is the following proposition.

\hfill

\proposition \label{_K-dR_inside_hK-dR_Dolbeault_Proposition_}
Let $M$ be a hyperk\"ahler manifold, $I, J, K$
the standard triple of induced complex structures,
$I\circ J = - J \circ I = K$, and 
$\g = \goth a \rtimes \c H$ the
hyperk\"ahler-de Rham superalgebra acting on $\Lambda^*(M)$.
Consider an 8-dimensional subspace $\g^{*,0}_I\subset \g$
spanned by the vectors
\begin{equation}\label{_g^*,0_I_basis_Equation_}
   L_\Omega,\Lambda_\Omega, H_\Omega, \6, \6^*, \6_J, \6_J^*, \Delta.
\end{equation}
 Then $\g^{*,0}_I$ is closed under the super-commutator.
Moreover, $\g^{*,0}_I$ is naturally
isomorphic to the K\"ahler-de Rham superalgebra.

\hfill

\remark 
In fact, the Dolbeault complex of a hyperk\"ahler manifold
is on most counts similar to the de Rham complex of a K\"ahler
manifold. This analogy was a driving engine behind
\cite{_Verbitsky:Hyperholo_bundles_}.

\hfill

{\bf Proof of \ref{_K-dR_inside_hK-dR_Dolbeault_Proposition_}:}
Let $W$ be the odd part of the hyperk\"ahler-de Rham superalgebra $\g$:
\[ W:= \langle d, d_I, d_J, d_K, d^*, d_I^*, d_J^*, d_K\rangle\]
Then, $W$ admits a Hodge decomposition
\[ W = W^{1,0} \oplus W^{0,1} \oplus W^{-1,0} \oplus W^{0,-1}
\]
with each component 2-dimensional, and generated by 
$\6$, $\6_J$, and their complex adjoint and Hermitian
conjugate as follows:
\begin{align*} 
{} &W^{1,0} = \langle \6, \6_J\rangle,  \ \ \  
W^{0,1} =\langle \bar\6, \bar\6_J\rangle,\\
{} & W^{-1,0} = \langle \6^*, \6_J^*\rangle,  \ \ \
W^{0,-1} = \langle \bar\6^*, \bar\6_J^*\rangle.
\end{align*}
Clearly, the odd part of $\g^{*,0}_I$ coinsides with
$W^{1,0}\oplus W^{-1,0}$. Since the $\goth{sl}(2)$-triple
$\langle L_\Omega,\Lambda_\Omega, H_\Omega \rangle$
maps $(p,q)$-forms to $(p',q)$-forms, it preserves
$W^{1,0}\oplus W^{-1,0}$. The pairwise 
super-commutators of odd
elements in $\g$ are proportional
to the Laplacian $\Delta$. Therefore, the space
$\g^{*,0}_I$ is closed under super-commutator.

To prove that $\g^{*,0}_I$ is isomorphic to the K\"ahler-de Rham superalgebra,
we use \ref{_Kah_deRha_formally_Theorem_}. It suffices to show
that the vectors $\6, \6_J$ anticommute\footnote{That is, commute in
the sense of superalgebra.} with themselves and
with each other, and commute with $L_\Omega$.

These vectors anticommute by Kodaira relations 
(see \ref{_hk-dR-acts_on_forms_Theorem_}). Since $\Omega$
is closed under the differentials $d, d_I, d_J, d_K$, we have
\[ [L_\Omega, \6]= [L_\Omega, \6_J]= 0. \]
Applying \ref{_Kah_deRha_formally_Theorem_},
we find that $\g^{*,0}_I$ is isomorphic to 
the K\"ahler-de Rham superalgebra. 
\ref{_K-dR_inside_hK-dR_Dolbeault_Proposition_}
is proven.
\endproof

\hfill

In HKT-geometry, $d\Omega=0$ and 
the hyperk\"ahler-de Rham superalgebra
does not act on $\Lambda^*(M)$. However, due to the HKT-relation
$\6\Omega=0$,
there exists the
K\"ahler-de Rham superalgebra action on $\Lambda^{*,0}(M,I)$. 
This structure is the main object of this paper.


\section{Lie superalgebra of an HKT-manifold}
\label{_superalh_HKT_Section_}


Let $M$ be a hypercomplex Hermitian manifold, $I$
an induced complex structure, and $\omega_I$ the
corresponding non-degenerate 2-form (\ref{_induced_c_str_Definition_}).
Denote by $L_I:\; \Lambda^*(M) \arrow \Lambda^{*+2}(M)$
the operator of multiplication by $\omega_I$, and let
$\Lambda_I$ the Hermitian adjoint operator
(Subsection 
\ref{_Dolbeault_hype_KdR_Subsection_}).  Clearly, 
the $\goth{so}(1,4)$-action 
(\ref{_SO(4,1)_Theorem_}) is valid 
for hypercomplex Hermitian manifolds
as well. However, there is no reason for 
Kodaira relations to hold; this is why 
there is no Hodge decomposition 
and no $\goth{so}(1,4)$-action
on the cohomology of a 
hyperk\"ahler manifold.

However, the HKT relation $\6\Omega=0$
(\ref{_HKT_via_diff_of_Omega_Theorem_} (i)) 
can be directly translated
to a statement about the Lie superalgebra.

Fix the standard quaternion basis $I, J, K$,
\[ I^2 = J^2 = K^2 =-1, I J = - JI = K, \]
Let $\Omega:= \frac 1 2 (\omega_J + \1\omega_K)$ be the corresponding
$(2,0)$-form (\ref{_induced_c_str_Definition_}),
\[ L_\Omega:\; \Lambda^*(M) \arrow \Lambda^{*+2}(M)\]
the operator of multiplication by $\Omega$,
and \[ \Lambda_\Omega:\; \Lambda^{*}(M) \arrow \Lambda^{*-2}(M)\] the 
Hermitian adjoint operator: $\Lambda_\Omega= *L_\Omega*$. 

From \ref{_SO(4,1)_Theorem_} it follows that the
commutator $H_\Omega:= [ L_\Omega, \Lambda_\Omega]$
is a scalar operator, mapping a $(p,q)$-form 
$\eta$ to $(n-p)\eta$, where $n = \dim_\H M$.
Therefore, $L_\Omega$, $\Lambda_\Omega$, $H_\Omega$
form an $\goth{sl}(2)$-triple on the standard Lie algebra
$\goth{so}(1,4)$ acting on $\Lambda^*(M)$.

\hfill

Consider the operator $J$ acting on differential 
forms as
\[ J(dx_1 \wedge dx_2 \wedge ... \wedge dx_n) =
   J(dx_1) \wedge J(dx_2) \wedge ... \wedge J(dx_n)
\]

where $dx_1\in \Lambda^1(M)$, and the action of $J$ on
$\Lambda^1(M)$ is determined from the quaternion structure.
Since $J$ and $I$ anticommute, $J$ maps $(p,q)$-forms
to $(q,p)$-forms.

Let $\6:\; \Lambda^{p,q}(M) \arrow \Lambda^{p+1,q}(M)$
be the Dolbeault differential, and 
\[ \6_J:\; \Lambda^{p,q}(M) \arrow \Lambda^{p+1,q}(M) \]
the differential obtained as
\[ \6_J := - J \circ \bar \6 \circ J. \]

\hfill

\proposition\label{_commu_L_Omega_6_J_Proposition_}
Let $M$ be a hypercomplex Hermitian manifold, 
with $L_\Omega$, $\6$ and $\6_J$ defined as above.
Then the HKT-relation $\6\Omega=0$
(\ref{_HKT_via_diff_of_Omega_Theorem_} (i)) is equivalent to
the following Lie algebra relations
\begin{equation}\label{_L_Omega_commu_w_6_Equation_} 
  [L_\Omega, \6] = [L_\Omega, \6_J] =0.
\end{equation}

{\bf Proof:} Clearly, $\6\Omega=0$ is equivalent to
$[L_\Omega, \6]=0$. However, $J(\Omega) = \bar \Omega$,
hence \[ \6_J\Omega = J(\bar\6\bar\Omega) = J(\overline{\6\Omega}).\]
Therefore, $\6_J\Omega=0$ if and only if $\6\Omega=0$.
We obtain that $\6\Omega=0$ implies $[L_\Omega, \6_J] =0$.
This proves \ref{_commu_L_Omega_6_J_Proposition_}.
\endproof

\hfill

\corollary \label{_K-dR_on_HKT_Corollary_}
Let $M$ be an HKT-manifold, 
with $L_\Omega$, $\Lambda_\Omega$, $H_\Omega$, 
$\6$ and $\6_J$ defined as above. Consider the
Lie superalgebra $\g$ generated by these operators.
Then $\g$ is isomorphic to the K\"ahler-de Rham superalgebra.

\hfill

{\bf Proof:} By \ref{_commu_L_Omega_6_J_Proposition_},
we have \[ [L_\Omega, \6] = [L_\Omega, \6_J] =0.\]
By \eqref{_6_6_J_commute_Equation_}, $\6$ and $\6_J$
(anti-)commute:
\[ \{\6, \6_J\} =0.\]
In Subsection \ref{_Dolbeault_hype_KdR_Subsection_}, 
we have seen that $H_\Omega$ acts on $(p, q)$-forms
as a muptiplication by $p-2n$. Therefore,
\[ [ H_\Omega, \6] = \6, \  \ [ H_\Omega, \6_J] = \6_J.
\]
From these relations and
\ref{_Kah_deRha_formally_Theorem_}
we obtain that $\g$ is isomorphic to the K\"ahler-de Rham superalgebra.
\endproof


\section{Kodaira relations in a differential graded algebra}
\label{_Kodaira_Section_}


\subsection{Differential operators in graded commutative rings}

Let $A$ be a graded commutative  algebra over $\C$. We define
the space $D_i(A)\subset \Hom_\C(A,A)$ of {\bf differential operators of order 
$\leq i$} recursively as follows.

\hfill

\definition\label{_diffe_ope_on_alge_Definition_}
(A. Grothendieck)
\begin{description}
\item[(i)] The space $D_0(A)$ of zero-order differential operators
is canonically identified with $A$ acting on 
itself by multiplication.
\item[(ii)] $D_i(A)$ is the space of all $\C$-linear maps
$D:\; A\arrow A$ such that the supercommutator $[a, D]$ belongs
to $D_{i-1}(A)$, for all $a\in D_0(A)$.
\end{description}

The differential operators form a Lie superalgebra.

\hfill

We have 
\[ 
   D_0(A) \subset D_1(A) \subset D_2(A) \subset ...
\]
The following lemma is well known

\hfill

\lemma\label{_commu_diffe_order_Lemma_}
We have
\[ [ D_i(A), D_j (A)] \subset D_{i+j-1}(A).
\]

\hfill

As usually, we denote the parity 
of an element $v$ of a graded vector
space by $\tilde v \in \{0,1\}$;
$\tilde v=0$ if $v$ is even, and
$\tilde v=1$ if $v$ is odd.

\hfill

The notion of differential operators can be extended to $A$-modules
in a usual fashion.

\hfill

\definition
Let $M$, $N$ be $A$-mudules. We define the space 
\[ D_i(M,N)\subset \Hom_\C(M,N)\]
of differential operators of order $\leq i$ recursively
as follows.

\begin{description}
\item[(i)] The space $D_0(M,N)$ of zero-order differential operators
is $\Hom_A(M,N)$ 
\item[(ii)] $D_i(M,N)$ is the space of all $\C$-linear maps
$D:\; M\arrow N$ such that for all $a\in A$, 
the supercommutator 
\[ [a, D](m) := a \cdot D(m) + (-1)^{\tilde a \tilde D} D(a\cdot m)
\]
 belongs
to $D_{i-1}(M,N)$, for all $a\in A$.
\end{description}

\hfill

\remark\label{_diffe_on_alge_and_usual_diffe_Remark_} 
Let $X$ be a smooth manifold, 
$A := \Lambda^*(X)$, $A_0:= C^\infty(X)$, and $M$, $N$ 
vector bundles on $X$, equipped with a structure
of $A$-module. Denote $M$, $N$, 
considered as $C^\infty(X)$-modules,
by $M_{C^\infty}$, $N_{C^\infty}$. Clearly, 
the $A_0$-differential operators $D_i(M_{C^\infty}, N_{C^\infty})$
are identified with the usual differential operators 
on the vector bundles corresponding to $M$, $N$. 
Since $A$-differential operator is necessarily
an $A_0$-differential operator, the space
$D_i(M,N)$ is a subspace of the space of 
$D_i(M_{C^\infty}, N_{C^\infty})$ of differential
operators on vector bundles. Generally speaking, 
the embedding
\[ D_i(M,N) \hookrightarrow D_i(M_{C^\infty}, N_{C^\infty})
\]
is proper.

\hfill

The main example we are working with is the following.

\hfill

\proposition \label{_d^*_seco_orde_Proposition_}
Let $X$ be a Riemannian manifold, $A = \Lambda^*(X)$
the algebra of differential forms on $X$, 
$*$ the Hodge star operator, and 
\[ d^*:= - * d *:\; A \arrow A
\]
the operator obtained as a Hermitian conjugate 
of the de Rham differential $d$.
Then $d^*$ is a second order differential
operator on the graded commutative algebra $A$:
$d^*\in D_2(A)$.

\hfill

{\bf Proof:} Consider 
the tensor product $A \otimes_{\C^\infty X} \Lambda^1(X)$
as an $A$-module. Let 
\[\nabla:\; A \arrow A \otimes_{\C^\infty X} \Lambda^1(X)
\] 
be the Levi-Civita connection operator. The Leibniz
formula implies
\[ \nabla\in D_1(A, A \otimes_{\C^\infty X} \Lambda^1(X)).
\]
Since $\nabla$ is torsion-free, we have 
$d^* = \nabla \circ \iota$, where
\[ \iota :\; A \otimes_{\C^\infty X} \Lambda^1(X) \arrow A 
\]
is the operator of ``inner multiplication'':
\[ \iota(\eta \wedge v) = - *((*\eta) \wedge v), 
\]
for all $\eta \in \Lambda^* (X)$, $v\in \Lambda^1 (X)$.

To prove \ref{_d^*_seco_orde_Proposition_},
it remains to show that the map
\[ \iota :\; A \otimes_{\C^\infty X} \Lambda^1(X) \arrow A 
\]
is a first-order differential operator of $A$-modules.
This is clear, because
\begin{multline*} \iota (dx_{1}\wedge dx_{2} \wedge ... \wedge dx_{k}\otimes v)
  \\ = \sum_{l=1}^k (-1)^{l-1} dx_{1}
    \wedge dx_{2} \wedge ... \wedge dx_{{l-1}} \wedge dx_{{l+1}}
    ... \wedge dx_{k} \cdot (v, dx_l)
\end{multline*}
where $v$ and $dx_i$ are $1$-forms and $(\cdot, \cdot)$ is
the Riemannian form.
\endproof

\hfill

\corollary \label{_commu_d^*_L_first_orde_Corollary_}
Let $\lambda$ be a differential form on $X$, 
$L_\lambda:\; A^* \arrow A^*$ the operation of
exterior  multiplication by $\lambda$, and 
$D:= [L_\lambda, d^*]\in D_*(A)$ the super-commutator
of $L_\lambda$ and $d^*$. Then $D\in D_1 (A)$.

\hfill

{\bf Proof:} Follows from \ref{_commu_diffe_order_Lemma_}
and \ref{_d^*_seco_orde_Proposition_}. \endproof

\hfill

\subsection{Kodaira relations and differential operators on 
hypercomplex Hermitian manifolds}

Let $M$ be a hypercomplex Hermitian manifold,
$I$, $J$, $K$ the standard triple of induced complex structures,
and $\Omega\in \Lambda^{2,0}_I(M)$.
the standard non-degenerate $(2,0)$-form on $(M,I)$. CVonsider
be the Dolbeault differential 
\[ \6:\; \Lambda^{p,0}_I(M) \arrow
   \Lambda^{p+1,0}_I(M)
\]
on $(M,I)$, and let
$ \6^* := -* \6 *$ be its Hermitian adjoint. Consider the
commutator
\[ \delta_J^*:= - [L_\Omega, \6^*].
\]

\hfill

\proposition \label{_commu_L_6^*_Proposition_}
In the above assumptions, the map
\[ \delta_J^*:\; \Lambda^{p,0}_I(M) \arrow \Lambda^{p+1,0}_I(M)
\]
is a first order differential operator 
(in the usual sense, that is, over $C^\infty M$)  
with the same symbol as the operator
\[ \6_J:\; \Lambda^{p,0}_I(M) \arrow \Lambda^{p+1,0}_I(M) \]
(Subsection \ref{_qD-bico_Subsection_}). Moreover, 
\[ \delta_J^*\in D_1(\Lambda^{*,0}_I(M)),\]
where $D_1(\Lambda^{*,0}(M))$ is the space
of first order differential operators
on the graded commutative algebra $\Lambda^{*,0}_I(M)$
(\ref{_diffe_ope_on_alge_Definition_}).\footnote{
The differential operators over $C^\infty M$ and
the differential
operators on the graded commutative algebra $\Lambda^{*,0}_I(M)$
are two distinct notions. See \ref{_diffe_on_alge_and_usual_diffe_Remark_} 
for a more explicit statement.}

\hfill

{\bf Proof:} A hypercomplex Hermitian manifolds admits
a first order approximation by a flat hyperk\"ahler
manifold $M' = {\Bbb H}^n$. On $M'$, we have
\[ [L_\Omega, \6^*] = -\6_J. \]
(\ref{_hk-dR-acts_on_forms_Theorem_}). 
The symbol of $\delta_J^*$
does not change if we replace $M$ by its first-order
approximation. This proves the first statement of
\ref{_commu_L_6^*_Proposition_}. The second statement
of \ref{_commu_L_6^*_Proposition_} is proven in exactly
the same way as \ref{_commu_d^*_L_first_orde_Corollary_}.
\endproof

\hfill

In assumptions of \ref{_commu_L_6^*_Proposition_},
consider the operator 
\[ h:= \delta_J^*- \6_J. \]
By \ref{_commu_L_6^*_Proposition_}, $h$ 
is a $\C^\infty M$-linear operator on $\Lambda^*(M)$,
which is a differential operator of first order over the
commutative graded algebra $\Lambda^{*,0}_I(M)$.
By definition of first order differential operators,
we have
\begin{multline}\label{_h_first_order_Equation_}
h(\eta_1\wedge \eta_2) = \\
(-1)^{\tilde \eta_2} \eta_1 \wedge h(\eta_2) 
+ h(\eta_1) \wedge \eta_2 + 
(-1)^{\tilde \eta_1+\tilde \eta+2}\eta_1\wedge \eta_2\wedge h(1).
\end{multline}
By \eqref{_h_first_order_Equation_}, $h$ is determined
by the values it takes on $\Lambda^0(M)$ and $\Lambda^{1,0}_I(M)$.
This allows to compute $h$ explicitly, as follows.

\hfill

Let $M$, $\dim_\H M =n$ be an HKT manifold,
 $I$, $J$, $K$ the standard quaternion triple,
$\Omega\in \Lambda^{2,0}_I(M)$ the corresponding
$(2,0)$-form, and 
\[ L_\Omega:\; \Lambda^{*,0}_I(M)\arrow 
   \Lambda^{*+2,0}_I(M)
\] 
the operator of multiplication by $\Omega$.
Let $\Omega^n$ be the corresponding 
nowhere degenerate section of the canonical
bundle. Since $\bar\6 \Omega^n$ is a $(2n,1)$-form,
and any $(2n,0)$-form on $(M,I)$ is proportional 
to $\Omega^n$, the form $\bar\6 \Omega^n$
is proportional to $\Omega^n$:
\[ \bar\6 \Omega^n = \bar \theta\wedge \Omega^n \]
(we denote the corresponding (0,1)-form by $\bar\theta$).
Let \[ J:\; \Lambda^{p,q}_I(M) \arrow \Lambda^{q,p}_I(M)\]
be the standard multiplicative map
acting on $\Lambda^1(M)$ as $J$. Denote
the form $J(\bar\theta) \in \Lambda^{1,0}_I(M)$ by $\theta_J$.

\hfill

\theorem\label{_delta_J_explicitly_Theorem_}
Let $M$ be an HKT-manifold.
In the above assumptions,
consider the operator
\[ \delta_J^*:= - [L_\Omega, \6^*].
\]
(\ref{_commu_L_6^*_Proposition_}). Let
\[ h:= \delta_J^*- \6_J.
\]
be the corresponding $C^\infty(M)$-linear operator. Then, for all
$\eta\in \Lambda^{*,0}_I(M)$, we have
\[h(\eta) =\theta_J\wedge  \eta ,\]
where $\theta_J$ is a 1-form, determined by the 
hypercomplex Hermitian structure as above.

\hfill

We prove \ref{_delta_J_explicitly_Theorem_}
in Subsection \ref{_*_on_HKT_and_6_Subsection_}.


\section{The Hodge $*$-operator on hypercomplex Hermitian manifolds}
\label{_*_Section_}


\subsection{The Hodge $*$-operator explicitly}

In this Subection, we perform explicit calculations related to the
Hodge $*$-operator on a hypercomplex Hermitian manifold. These
calculations are purely linear-algebraic; in fact, one can
compute everything on a quaternionic Hermitian space.

\claim\label{_*_expli_Claim_}
Let $M$ be a hypercomplex Hermitian manifold,
$I$, $J$, $K$ the standard quaternion triple, and
$\Omega\in \Lambda^{2,0}_I(M)$ the corresponding
$(2,0)$-form on $(M,I)$. Then
\begin{description}
\item[(i)] $*1 = (\frac{1}{n!})^2 \Omega^n \wedge \bar\Omega^n$, \ \  $n=\dim_\H M$
\item[(ii)] $*\Omega = 2n(\frac{1}{n!})^2\Omega^{n-1}\wedge  \bar\Omega^n$
\item[(iii)] For any $(1,0)$-form $\eta$, we have
\[ *\eta = - (\frac{1}{n!})^2 \Omega^{n-1} \wedge \bar\Omega^n \wedge J(\bar \eta),
\]
where 
\[ J:\; \Lambda^{p,q}_I(M) \arrow \Lambda^{q,p}_I(M)\]
is a standard multiplicative operator on differential forms
acting on $\Lambda^1(M)$ as $J$.
\end{description}
{\bf Proof:} \ref{_*_expli_Claim_} (i) is clear, 
because $*1$ is by definition equal to the Riemannian 
volume form of $M$, and
\[ \Vol(M) = \frac{\Omega^n \bar\Omega^n}{{n!}^2}.
\]
To prove \ref{_*_expli_Claim_} (ii),
we write 
\begin{equation}\label{_*_Omega_via_Lambda_Equation_}
   * \Omega = *(L_\Omega 1) = \Lambda_\Omega (* 1) 
   = \Lambda_\Omega L_\Omega^n \bar\Omega^n \cdot \frac{1}{(n!)^2}
\end{equation}
(the last equation follows from \ref{_*_expli_Claim_} (i)).
Since $L_\Omega$, $\Lambda_\Omega$, $H_\Omega$ form
an $\goth{sl}(2)$-triple, 
and $\bar\Omega^n $ is a lowest vector of a weight $n$ representation,
we have 
\begin{equation}\label{_Omega_commu_lowest_we_Equation_}
  \Lambda_\Omega L_\Omega^n \bar\Omega^n = 2n\Omega^{n-1}\bar\Omega^n.
\end{equation}

Comparing \eqref{_Omega_commu_lowest_we_Equation_}
and \eqref{_*_Omega_via_Lambda_Equation_}, we obtain
\ref{_*_expli_Claim_} (ii).

\hfill

To prove \ref{_*_expli_Claim_} (iii), we use the following
elementary lemma.

\hfill

\lemma \label{_[Lambda_eta_L_Omega_](commu)_Lemma_}
Let $M$ be a hypercomplex Hermitian manifold,
$I$, $J$, $K$ the standard quaternion triple, and
$\Omega\in \Lambda^{2,0}_I(M)$ the corresponding
$(2,0)$-form on $(M,I)$. Consider an arbitrary
(1,0)-form $\eta \in \Lambda^{1,0}_I(M)$.
Denote by 
\[ L_\eta:\; \Lambda^{*,0}_I(M)\arrow \Lambda^{*+1,0}_I(M)
\] the operator of exterior
multiplication by $\eta$, and let
$\Lambda_\eta:= -* L_\eta *$ be the Hermitian
adjoint operator (so-called ``inner multiplication by $\eta$'').
Let $L_\Omega$ be the  operator of exterior
multiplication by $\Omega$. Then
\[ 
  [ L_\Omega, \Lambda_\eta] = L_{J(\bar \eta)}
\]
where 
\[ J:\; \Lambda^{p,q}_I(M) \arrow \Lambda^{q,p}_I(M)\]
is the standard operator defined above.

\hfill

{\bf Proof:}
Proven by a local computation. \endproof

\hfill

Return to the proof of \ref{_*_expli_Claim_} (iii).
By \ref{_[Lambda_eta_L_Omega_](commu)_Lemma_} and
\ref{_*_expli_Claim_} (i), we have 
\begin{align*} *(\eta) {} &= * L_\eta 1 = - \Lambda_\eta (*1) 
   = \frac{1}{(n!)^2} \Lambda_\eta L_\Omega^n \bar \Omega^n \\{} &=
   -\frac{1}{(n!)^2} L_{J(\bar \eta)}L_\Omega^{n-1} \bar \Omega^n =
    -\frac{1}{(n!)^2} J(\bar \eta) \wedge \Omega^{n-1}\wedge
\Omega^n.
\end{align*}
This is exactly the statement of \ref{_*_expli_Claim_} (iii).
We proved \ref{_*_expli_Claim_}. \endproof

\subsection{The Hodge $*$-operator on HKT-manifolds and
Dolbeault differential}
\label{_*_on_HKT_and_6_Subsection_}

We work in assumptions of \ref{_delta_J_explicitly_Theorem_}.
By \ref{_commu_L_6^*_Proposition_}
(see also \eqref{_h_first_order_Equation_}), to prove 
\ref{_delta_J_explicitly_Theorem_}, it suffices to check
\[ 
\delta^*_J\eta - \6_J \eta = J(\bar \theta)\wedge \eta,
\]
for $\eta\in \Lambda_I^{*,0}(M)$ a $0$-form or an $(1,0)$-form.
The operator $h= \delta^*_J-\6_J$ is $C^\infty$-linear.
Since $\6_J^2=0$, the $\6_J$-closed forms generate
the space of all forms over $C^\infty(M)$.
Therefore, to prove \ref{_delta_J_explicitly_Theorem_}, 
it suffices to check
\begin{equation}
h(\eta) = J(\bar \theta)\wedge \eta,
\end{equation}
for $\6_J(\eta) =0$, $\eta \in \Lambda_I^{i,0}(M)$, $i=0,1$.
For $\eta$ a $\6_J$-closed function,
we have
\[ \delta^*_J\eta = - [L_\Omega, \6^*] \eta = \6^* (\eta\Omega).
\]
($\6^*\eta=0$ because $\eta$ is 0-form).
Using \ref{_*_expli_Claim_} (ii), we obtain
\begin{equation}\label{_delta_of_fu_expli_Equation_}
   \delta^*_J\eta = 
   -2n \frac{1}{(n!)^2} * \6 (\Omega^{n-1}\wedge \bar\Omega^n).
\end{equation}
Using the definition of $\theta$ and the HKT equation $\6 \Omega=0$,
we express the right hand side of  
\eqref{_delta_of_fu_expli_Equation_} 
through $\theta$ as follows:
\[ 
\delta^*_J\eta = 
-2n \frac{1}{(n!)^2}\eta *(\theta \wedge \Omega^{n-1}\wedge \bar\Omega^n).
\] 
We obtain 
\[ h(\eta) = 
   -\Lambda_\theta \bigg( *\bigg[2n \frac{1}{(n!)^2}\Omega^{n-1}\wedge \bar\Omega^n\bigg]\bigg)\eta
\]
where $\eta$ is an arbitrary 0-form. The expression in brackets
is equal to $*\Omega$ as follows from 
\ref{_*_expli_Claim_} (ii). Therefore,
\begin{equation}\label{_h_func_Equation_}
  h(\eta) = -\Lambda_\theta(\Omega)\cdot \eta.
\end{equation}
for any $\eta\in \Lambda^0(M)$.
By  \ref{_[Lambda_eta_L_Omega_](commu)_Lemma_},
\begin{equation}\label{_commu_Lambda_theta_and_L_Ome(1)_Equation_}
   -\Lambda_\theta(\Omega) = [L_\Omega, \Lambda_\theta] (1) 
   = L_{J(\bar\theta)}(1) = J(\bar\theta).
\end{equation}
Comparing \eqref{_h_func_Equation_} and 
\eqref{_commu_Lambda_theta_and_L_Ome(1)_Equation_},
we obtain
\[ h(\eta) = J(\bar\theta)\cdot \eta,
\]
where $\eta$ is an arbitrary 0-form.
Now assume that
$\eta\in \Lambda^{1,0}_I(M)$, $\6_J(\eta)=0$.
By \ref{_*_expli_Claim_} (iii), we have
\begin{equation} \label{_delta^*_J_via_commu_6_Lambda_Equation_}
\begin{aligned}
  \delta^*_J \eta {} & = -* [\Lambda_\Omega, \6] *\eta\\ {} &=
   -\frac{1}{(n!)^2} *\bigg( [\Lambda_\Omega, \6]\bigg( 
          \Omega^{n-1}\wedge \bar\Omega^n\wedge J(\bar\eta)\bigg)\bigg).
\end{aligned}
\end{equation}
Since $\eta$ is $\6_J$-closed, we have 
$\6(J(\bar\eta))=0$. 
We have $\6\Omega=0$ by HKT-equations.
Since $J(\bar\eta)\wedge \bar\Omega^n$ is a 
lowest weight vector with respect to the
$\goth{sl}(2)$-action generated by $L_\Omega$, $\Lambda_\Omega$,
$H_\Omega$, the form
\[ \Lambda_\Omega\left(\Omega^{n-1}\wedge 
   \bar\Omega^n\wedge J(\bar\eta)\right)
\]
is proportional to
\[
\Omega^{n-2}\wedge 
   \bar\Omega^n\wedge J(\bar\eta)
\]
Therefore, $\6$ acts on 
$\Omega^{n-1}\wedge \bar\Omega^n\wedge J(\bar\eta)$ 
and on $\Lambda_\Omega\left(\Omega^{n-1}\wedge 
   \bar\Omega^n\wedge J(\bar\eta)\right)$
as a multiplication by $\theta$, and 
we can replace $\6$ in 
\eqref{_delta^*_J_via_commu_6_Lambda_Equation_}
by $L_\theta$. We obtain
\[ \delta^*_J \eta =\frac{1}{(n!)^2} *\bigg( [\Lambda_\Omega, L_\theta] 
          \bigg(\Omega^{n-1}\wedge \bar\Omega^n\wedge J(\bar\eta)\bigg)\bigg) =
   [ L_\Omega, \Lambda_\theta] \eta.
\]
On the other hand, $[ L_\Omega, \Lambda_\theta]= L_{J(\bar\theta)}$
(\ref{_[Lambda_eta_L_Omega_](commu)_Lemma_}).
We obtain
\[ \delta^*_J \eta=J(\bar\theta)\wedge\eta,
\]
for any $\eta\in \Lambda^{1,0}_I(M)$, $\6_J(\eta)=0$.
This finishes the proof of 
\ref{_delta_J_explicitly_Theorem_}.
\endproof


\section{Hodge theory on $\Lambda^{p,0}(M)\otimes K^{1/2}$}
\label{_K^1/2_Section_}


\subsection{The normalized K\"ahler-de Rham superalgebra 
of an HKT manifold}
\label{_normalized_g_defi_Subsection_}

Let $M$ be an HKT manifold and 
\[ \goth g= \langle L_\Omega,
   \Lambda_\Omega, H_\Omega, \6, \6_J, \delta, \delta_J,
  \Delta\rangle
\]
be the corresponding K\"ahler-de Rham superalgebra 
(\ref{_K-dR_on_HKT_Corollary_}).
Conjugating each generator of $\goth g$ by $*$, we obtain
another superalgebra
\[ 
\goth g^*= \langle L_\Omega,
   \Lambda_\Omega, H_\Omega, \6^*, \6_J^*, \delta^*, \delta_J^*,
  \Delta^*\rangle
\]
which is naturally isomorphic to $\goth g$. 
However, these algebras are distinct; to work
with Hodge theory, we need to relate the 
Lie superalgebras and the Hodge $*$-operator.

In this section, we construct and study the ``normalized''
K\"ahler-de Rham superalgebra ${}^n\goth g$, which acts on
$\Lambda_I^{*,0}(M)$ in the same way as 
$\goth g$ and $\goth g^*$, and is fixed
by $*$.

Consider the 4-dimensional vector space ${}^nV$ 
spanned by the vectors
\begin{align*} 
{}^n \6 := \frac{\6+\delta^*}{2}, {} &\ \  {}^n \6_J := \frac{\6_J+\delta^*_J}{2},
\\ {}^n \6^* := \frac{\6^*+\delta^*}{2}, {} &\ \ 
  {}^n \6_J^* := \frac{\6_J^*+\delta_J}{2}\in \End (\Lambda_I^{*,0}(M)).
\end{align*}
Consider the standard $\goth{sl}(2)$-action on $\Lambda_I^{*,0}(M)$
associated with the $\goth{sl}(2)$-triple $L_\Omega,
   \Lambda_\Omega, H_\Omega$. Since this $\goth{sl}(2)$-action
exchanges $\6$ with $\delta_J$ and $\delta^*$ with $\delta^*_J$,
this action exchanges ${}^n \6$ with ${}^n \6_J^*$:
\[ [L_\Omega, {}^n \6_J^*] = {}^n \6, \ \ 
   [\Lambda_\Omega, {}^n \6_J^*]  = -{}^n \6.
\]
In other words, the standard 
$\goth{sl}(2)$-triple 
\[ \langle L_\Omega,
   \Lambda_\Omega, H_\Omega\rangle\] 
preserves the
2-dimensional space $\langle {}^n \6, {}^n \6_J^*\rangle$
spanned by ${}^n \6$, ${}^n \6_J^*$. Similarly,
it preserves the space 
$\langle {}^n \6^*, {}^n \6_J\rangle$.

We obtain that this $\goth{sl}(2)$-action 
preserves the space ${}^n V\subset\End (\Lambda_I^{*,0}(M))$,
spanned by these two 2-dimensional spaces defined above.

Let $\theta\in \Lambda_I^{1,0}(M)$ be a 1-form 
defined by the formula
\[ \6 \bar K = \theta \wedge \bar K,\]
where $\bar K\in \Lambda_I^{0, 2n}(M)$, $\bar K = \bar \Omega^n$,
$n = \dim_\H M$ the natural section of the
line bundle $\Lambda_I^{0, 2n}(M)$ determined by 
the standard nowhere degenerate $(0,2)$-form
$\bar \Omega$. Denote by $\theta_J$ the form
$J(\bar \theta)\in \Lambda_I^{1,0}(M)$.
By \ref{_delta_J_explicitly_Theorem_},
we have
\[ 
  {}^n\6 - \delta^* = \theta, \ \ \  {}^n\6_J - \delta^*_J = \theta_J.
\]
This implies 
\begin{equation}\label{_^n6-explicitly_Equation_}
{}^n\6 = \6 + \frac{1}{2} \theta, \ \ \ {}^n\6_J = \6_J + \frac{1}{2} \theta_J.
\end{equation}
By definition of $\theta$, we have
\begin{align*} 
  0 = \6 \6 \bar K {} & = \6(\theta\wedge \bar K)\\ {} & = 
 \6(\theta) \wedge\bar K - \theta\wedge \6 \bar K =
  \6(\theta) \wedge\bar K -\theta\wedge\theta\wedge\bar K \\{} & = 
  \6(\theta) \wedge\bar K.
\end{align*}
Therefore, $\theta$ is $\6$-closed, and we have 
\[ 
{}^n\6 ^2 = (\6 - \frac{1}{2} \theta)^2 =0.
\]
Twisting this equation with $J$, we obtain ${}^n\6_J ^2=0$.
Finally, $\6$ and $\6_J$ anticommute \eqref{_6_6_J_commute_Equation_}, 
and we have
\begin{multline}\label{_6_6_J_to_bar_K_Equation_}
   0 = \6\6_J \bar K + \6_J \6\bar K = \6(\theta_J \wedge \bar K)+
       \6_J(\theta \wedge\bar K) \\= 
       \6(\theta_J) \wedge \bar K+\6_J(\theta)\wedge \bar K
       -\theta_J \wedge\6\bar K - \theta \wedge\6_J\bar K.
\end{multline}
The last two terms of \eqref{_6_6_J_to_bar_K_Equation_}
are equal to $-\theta_J \wedge\theta\wedge\bar K$ and
$-\theta \wedge\theta_J\wedge\bar K$, therefore they
cancel each other. We reduced \eqref{_6_6_J_to_bar_K_Equation_}
to the equation
\[ 0 = \6(\theta_J) \wedge \bar K+\6_J(\theta)\wedge \bar K.\]
This implies
\[ \6\theta_J +\6_J\theta=0.\]
Using \eqref{_^n6-explicitly_Equation_},
we immediately obtain the anticommutation relation
\[ \{{}^n\6, {}^n\6_J\} =0.
\]
We have checked all the conditions of \ref{_Kah_deRha_formally_Theorem_}.
From \ref{_Kah_deRha_formally_Theorem_}, 
we obtain that the $\goth{sl}(2)$-triple 
 $\langle L_\Omega,  \Lambda_\Omega, H_\Omega\rangle$
and the 4-dimensional space 
\begin{multline*}
{}^nV:= 
\bigg \langle{}^n \6 := \frac{\6+\delta^*}{2}, \ \  {}^n \6_J := \frac{\6_J+\delta^*_J}{2},
\\ {}^n \6^* := \frac{\6^*+\delta^*}{2}, \ \ 
  {}^n \6_J^* := \frac{\6_J^*+\delta_J}{2}\bigg\rangle
  \subset \End (\Lambda_I^{*,0}(M)).
\end{multline*}
generate a K\"ahler-de Rham 
Lie superalgebra, denoted by ${}^n\goth g$.

\hfill

\definition\label{_normali_g_Definition_}
This Lie superalgebra called {\bf the normalized 
K\"ah\-ler-\-de Rham superalgebra of an HKT-manifold}.

\subsection[The normalized 
K\"ahler-de Rham su\-per\-algebra and Lefschetz theorem]{The normalized 
K\"ahler-de Rham su\-per\-algebra \\ and Lefschetz theorem}
\label{_Lefschetz_Subsection_}

The normalized 
K\"ahler-de Rham superalgebra ${}^n\goth g$
has the following geometrical interpretation. 

Let $M$ be an HKT-manifold, $\dim _\H M=n$, 
$I$, $J$, $K$ the standard triple of quaternions,
 ${\sf K}=\Lambda_I^{2n,0}(M)$ the canonical bundle of 
$(M, I)$ and $\goth V= \Omega^n$ its nowhere degenerate
section provided by the canonical nowhere degenerate
(2,0)-form $\Omega\in \Lambda_I^{2,0}(M)$. 
Since an HKT-manifold is Hermitian, the bundle
${\sf K}$  is equipped with a 
Hermitian metrics.
Consider the standard Hermitian connection $\nabla$
on ${\sf K}$ associated with this
metrics. Since ${\goth V}$ has constant length,
we have
\begin{equation}\label{_conne_via_K_Equation_}
\nabla({\goth V}) = \bar\6 ({\goth V}) + \overline{\bar\6 {\goth V}}, 
\end{equation}
where $\bar\6 :\; {\sf K} \arrow {\sf K}\otimes \Lambda_I^{0,1}(M)$
is the holomorphic structure on ${\sf K}$. Interpreting
${\goth V}$ as a $(2n,0)$-form on $M$, we find
\begin{equation}\label{_connec_via_theta_Equation_}
\bar\6 {\goth V} = \bar\theta\otimes {\goth V}, 
\end{equation}
where $\theta$
is the canonical $(0,1)$-form defined 
in Subsection \ref{_normalized_g_defi_Subsection_} 
(see also \ref{_delta_J_explicitly_Theorem_}).
Since ${\goth V}$ is a nowhere degenerate section of
${\sf K}$, ${\goth V}$ provides a $C^\infty$-trivialization
of this bundle. Let $\nabla_0$ be a flat
connection associated with this trivialization.
Comparing \eqref{_conne_via_K_Equation_} and
\eqref{_connec_via_theta_Equation_}, we find that $\nabla$
can be expressed via $\nabla_0$ as follows:
\[
\nabla= \nabla_0 +\theta+\bar\theta.
\]
In other words, $\theta$ is the $(1,0)$-connection form of 
the canonical bundle ${\sf K}$
associated with the $C^\infty$-trivialization 
provided by the $C^\infty$-section ${\goth V}\in {\sf K}$.

Let ${\sf K}^{1/2}$ be a square root of ${\sf K}$ 
determined by the above trivialization.
One can define ${\sf K}^{1/2}$ as a trivial
$C^\infty$-bundle with a holomorphic
structure defined by a connection
\[
\nabla_{1/2} = \nabla_0 +\frac{1}{2}\theta+\frac{1}{2}\bar\theta.
\]
Consider the (anti-)Dolbeault complex of ${\sf K}^{1/2}$:
\begin{equation}\label{_Dolb_in_K^1/2_Equation_}
{\sf K}^{1/2}
\stackrel{\nabla_{1/2}^{1,0}}\arrow {\sf K}^{1/2}\otimes \Lambda^{1,0}_I(M)
\stackrel{\nabla_{1/2}^{1,0}}\arrow {\sf K}^{1/2}\otimes \Lambda^{2,0}_I(M) ...
\end{equation}
where $\nabla_{1/2}^{1,0}$ is the $(1,0)$-component of $\nabla$.
Using the standard trivialization of ${\sf K}^{1/2}$,
we may identify ${\sf K}^{1/2}\otimes \Lambda^{*,0}_I(M)$
and $\Lambda^{*,0}_I(M)$.
By definition of $\nabla_{1/2}$, we have
\[ 
   \nabla_{1/2}^{1,0} = \6 + \frac{1}{2}\theta
\]
In other words, the twisted Dolbeault differential
\eqref{_Dolb_in_K^1/2_Equation_} is equal to the
normalized HKT differential ${}^n\6$ of \ref{_normali_g_Definition_}.

\hfill

Consider the action of 
$\goth{sl}(2)= \langle L_\Omega,  \Lambda_\Omega, H_\Omega\rangle$
on the complex \eqref{_Dolb_in_K^1/2_Equation_}.
Since $\goth{sl}(2)$, ${}^n\6$ and ${}^n\6^*$ are elements
of K\"ahler-de Rham superalgebra, they satisfy the
conditions of \ref{_Kah_deRha_formally_Theorem_}. In particular, 
the Laplace operator
\[ {}^n\Delta:= {}^n\6{}^n\6^* +{}^n\6^*{}^n\6
\]
commutes with the $\goth{sl}(2)$-action. On the other
hand, the cohomology of \eqref{_Dolb_in_K^1/2_Equation_}
are identified with the kernel of the Laplacian
${}^n\Delta$.
This proves the Lefschetz theorem for
the cohomology of the complex \eqref{_Dolb_in_K^1/2_Equation_}:

\hfill

\theorem\label{_Lefschetz_Theorem_}
Let $M$ be a compact HKT-manifold, and $I, J, K$ the standard 
triple of quaternions, and ${\sf K}^{1/2}$
be the square root of the canonical class of $(M,I)$
constructed as above. Consider its cohomology space 
$H^*(M, {\sf K}^{1/2})$.\footnote{This 
is the cohomology of the complex \eqref{_Dolb_in_K^1/2_Equation_}.}
Then 
\begin{description}
\item[(i)] Using the standard trivialization of ${\sf K}^{1/2}$,
we can identify the Dolbeault complex of the bundle ${\sf K}^{1/2}$ with 
the complex 
\begin{equation}\label{_^n6_complex_Equation_}
0 \arrow \Lambda_I^{0,0}(M) \stackrel{{}^n\6} \arrow 
\Lambda_I^{1,0}(M) \stackrel{{}^n\6} \arrow \Lambda_I^{2,0}(M) \stackrel{{}^n\6} \arrow ...
\end{equation}

\item[(ii)] The cohomology of ${\sf K}^{1/2}$, 
or, what is the same, the cohomology of 
the complex \eqref{_^n6_complex_Equation_}, are identified with
the kernel of the corresponding Laplacian
\[ {}^n \Delta = \{{}^n\6, \6^*\} \]

\item[(iii)] The $\goth{sl}(2)$-triple $\langle L_\Omega$,
$\Lambda_\Omega$, $H_\Omega\rangle$ commutes with the action of  
the normalized Laplacian ${}^n \Delta$. 

\item[(iv)]  This $\goth{sl}(2)$-action provides 
the ``Hard Lefschetz'' isomorphism
\[ L_\Omega^{n-i}:\; H^i({\sf K}^{1/2}) 
   \arrow H^{2n-i}({\sf K}^{1/2}).
\]
Together with the Serre's duality 
\[ H^{2n-i}({\sf K}^{1/2})\cong H^{i}({\sf K}^{1/2})^*, \]
this gives a canonical isomorphism
\begin{equation}\label{_pairing_on_coho_of_sq_of_K_Equation_} 
  H^i({\sf K}^{1/2}) \cong H^{i}({\sf K}^{1/2})^*.
\end{equation}

\item[(v)] Consider the map 
\[ J_c:\; \Lambda_I^{p,0}(M) \arrow\Lambda_I^{p,0}(M),\]
$J_c(\eta) = J(\bar \eta)$, where 
$J:\; \Lambda_I^{p,0}(M) \arrow \Lambda_I^{0,p}(M)$
is the standard multiplicative
map of differential forms 
associated with the induced complex structure $J$. 
Then $J_c$ commutes with the normalized Laplacian 
${}^n \Delta$. 

\item[(vi)] The pairing
\eqref{_pairing_on_coho_of_sq_of_K_Equation_} 
can be obtained explicitly as follows. Take $\eta$, 
$\eta'\in H^i({\sf K}^{1/2})$ ($i\leq n$).
Consider the cohomology class 
\[ \eta \wedge J_c(\eta') \in H^{2i}(M, {\sf K}),\] and let
\[ \alpha:=L_\Omega^{n-i}(\eta \wedge J_c(\eta') )\] be the
corresponding element in $H^{2n}(M, {\sf K})\cong \C$.
Then $\langle \eta, \eta'\rangle=\alpha$, where
$\langle \cdot,\cdot\rangle$ is the pairing 
\eqref{_pairing_on_coho_of_sq_of_K_Equation_}.

\end{description}

\endproof

\subsection{Harmonic spinors on HKT manifolds}

We work in assumptions and notations of 
\ref{_Lefschetz_Theorem_}. The canonical class
of a hypercomplex Hermitian manifold is topologically trivial:
take, for instance, a trivialization associated
with the section $\goth V$ of the canonical class
(Subsection \ref{_Lefschetz_Subsection_}).
Therefore, a hypercomplex Hermitian manifold 
admits a natural spinor structure. 

Let $\cal S$ be the spinor bundle of $M$. 
Then $\cal S$ is isomorphic to 
$\Lambda^{*,0}_I(M)\otimes {\sf K}^{1/2}$,
where ${\sf K}^{1/2}$ is the square root 
of the canonical class constructed above. 

It is well known that the Dirac operator
corresponds to ${}^n\6 + {}^n\6^*$ under the
identification
\[ \Lambda^{*,0}_I(M)\otimes {\sf K}^{1/2} \cong \cal S
\]
where ${}^n\6$ denotes the Dolbeault differential
\eqref{_^n6_complex_Equation_}. Therefore,
$H^*({\sf K}^{1/2})$ is naturally 
identified with a space
of harmonic spinors. Since the harmonic
spinors depend from a metric only
and not from the choice of a complex 
structure, we obtain the following
corollary.

\hfill

\corollary
In assumptions of \ref{_Lefschetz_Theorem_},
consider the space $H^*({\sf K}^{1/2})$.
Then $H^*({\sf K}^{1/2})$ is canonically
isomorphic to the space of harmonic spinors. 
Moreover, $H^*({\sf K}^{1/2})$ does
not depend from the choice of the basis
$I, J, K$ in quaternions.

\endproof

\hfill

\remark
Let $G$ be a compact semisimple Lie group
from the list \eqref{_Joyce_list_Equation_},
equipped with an HKT-structure as in
\ref{_Joyce_exa_Example_}. Then 
$H^*({\sf K}^{1/2})=0$. The algebraic-geometrical
computation needed for this result was performed 
by D. Kaledin.

\hfill

\hfill

{\bf Acknowledgements:}
This paper appeared as a result of fruitful
talks with D. Kaledin. I am grateful to
Y.-S. Poon, who told me of HKT-manifolds; 
 A. Losev, who suggested the close study of
Lie superalgebras arising from the natural DG-algebras
on manifolds; and M. Kontsevich,
who told me of M. Reid's conjecture. 
Many thanks  to S. Merkulov for his
interest and encouragement and R. Bielawski for 
interesting discussions. My gratitude to F. Bogomolov
and D. Kaledin, who explained me the Bott-Samelson
construction. G. Papadopoulos kindly send
some remarks on the history of the subject.

\hfill


\end{document}